\documentclass[12pt]{article}

\usepackage[latin1]{inputenc}

\usepackage{pgf,tikz}
\usetikzlibrary{arrows}
\usetikzlibrary{patterns}

\usepackage{amssymb}

\usepackage{amsmath,epsfig}
\usepackage{multirow} 
\usepackage{subfigure}

\newtheorem{Theorem}{Theorem}[section]
\newtheorem{Corollary}[Theorem]{Corollary}

\newtheorem{Definition}[Theorem]{Definition}
\newtheorem{Lemma}[Theorem]{Lemma}
\numberwithin{equation}{section}

\newcommand{\C}{{\mathbb C}}

\newcommand{\beq}{\begin{equation}}
\newcommand{\eeq}{\end{equation}}

\def\qed{\hfill$\Box$\\ \medskip}

\def\mybox{\hbox to 12.0pt}

\def\mybigbox{\hbox to 35.0pt}
\def\myverybigbox{\hbox to 60.0pt}

\def\g{g}
\def\gl{}
\def\gl{}
\def\sp{sp}
\def\so{so}
\def\oo{so}

\def\ov{\overline}

\def\wgt{{\rm wgt}}

\def\sgn{{\rm sgn}}

\def\ch{{\rm ch}\,}
\def\u{{\bf u}}
\def\v{{\bf v}}
\def\w{{\bf w}}
\def\x{{\bf x}}
\def\y{{\bf y}}
\def\z{{\bf z}}

\def\u{{\bf u}}
\def\a{{\bf a}}

\def\0{{\bf 0}}
\def\ovx{{\bf\ov{\x}}}
\def\ovy{{\bf\ov{\y}}}

\def\1{{\bf 1}}

\def\+{\!\!+\!\!}

\newcommand{\YT}[4]{
\vcenter{\hbox{
\begin{tikzpicture}[x={(0in,-#1)},y={(#2,0in)}] 
\foreach \k [count=\i] in {#4} {
 \foreach \e [count=\j] in \k {
  \draw (\i,\j) rectangle +(-1,-1);
  \draw (\i-0.5,\j-0.5) node {$#3\e$};
 }
}
\end{tikzpicture}
}}
}

\newcommand{\wideYT}[4]{
\vcenter{\hbox{
\begin{tikzpicture}[x={(0in,-#1)},y={(#2,0in)}] 
\foreach \k [count=\i] in {#4} {
 \foreach \e [count=\j] in \k {
  \draw (\i,\j) rectangle +(-1,-1);
  \draw (\i-0.5,\j-0.5) node {$#3\e$};
 }
}
\end{tikzpicture}
}}
}

\newcommand{\SYT}[4]{
\vcenter{\hbox{
\begin{tikzpicture}[x={(0in,-#1)},y={(#2,0in)}] 
\foreach \k [count=\i] in {#4} {
 \foreach \e [count=\j] in \k {
  \draw (\i,\j+\i-1) rectangle +(-1,-1);
  \draw (\i-0.5,\j+\i-1-0.5) node {$#3\e$};
 }
}
\end{tikzpicture}
}}
}

\usepackage{color}
\def\red{\textcolor{red}}
\def\blue{\textcolor{blue}}
\def\magenta{\textcolor{magenta}}

\title{Factorial characters of some classical Lie groups}

\author{
Ang\`ele M. Hamel\thanks{ 
Department of Physics and Computer Science,
Wilfrid Laurier University,
Waterloo, Ontario, N2L 3C5, Canada ({\tt ahamel@wlu.ca})}
\and 
Ronald C. King\thanks{
Mathematical Sciences, University of Southampton, 
Southampton SO17 1BJ, England ({\tt R.C.King@soton.ac.uk})}}

\begin{document}

\maketitle

\begin{abstract}
A definition is offered of the factorial characters of the general linear group, the symplectic group
and the orthogonal group in an odd dimensional space. It is shown that these characters satisfy
certain flagged Jacobi-Trudi identities. These identities are then used to give combinatorial
expressions for the factorial characters: first in terms of a lattice path model and then in terms of
the well known tableaux associated with the classical groups. Factorial $Q$-functions are then
defined in terms of three sets of primed shifted tableaux, and shown to satisfy Tokuyama type identities 
in each case.

\end{abstract}

\section{Definition of factorial characters}\label{sec-fact-char}

A class of symmetric polynomials $t_\lambda(\x)$ labelled by partitions $\lambda=(\lambda_1,\lambda_2,\ldots,\lambda_n)$
was introduced by Biedenharn and Louck~\cite{BL89} as a new integral basis of the ring of all symmetric polynomials in the 
parameters $\x=(x_1,x_2,\ldots,x_n)$. These new symmetric but inhomogeneous polynomials were studied further
and called factorial Schur functions by Chen and Louck~\cite{CL93}. They were given a more general form by 
Goulden and Greene~\cite{GG94} and by Macdonald~\cite{Mac92} who introduced the notation $s_\lambda(\x\,|\,\a)$,
where $\a=(a_1,a_2,\ldots)$ is an infinite sequence of factorial parameters. Macdonald gave a definition of
factorial Schur functions $s_\lambda(\x\,|\,\a)$ as a ratio of determinants, exactly analagous to that applying
to ordinary Schur functions $s_\lambda(\x)$. Since Schur functions, through 
their determinantal definition, can be identified with characters of $GL(n,\C)$, it is not unreasonable to
refer to factorial Schur functions as factorial characters of $GL(n,\C)$. What we then aim to offer here is
a definition of factorial characters of the classical groups $Sp(2n,\C)$ and $SO(2n+1,\C)$. We defer
the more difficult case of $SO(2n,\C)$ for consideration elsewhere.

For each of the three groups $G=GL(n,\C)$, $Sp(2n,\C)$ and $SO(2n+1,\C)$ there exists a finite dimensional 
irreducible representation $V_G^\lambda$ of highest weight $\lambda$, where $\lambda=(\lambda_1,\lambda_2,\ldots,\lambda_n)$ 
is a partition of length $\ell(\lambda)\leq n$. Its character may be denoted by $\ch V_G^\lambda(\z)$ where $\z=(z_1,z_2,\ldots,z_N)$
is a suitable parametrisation of the $N$ eigenvalues of the group elements of $G=GL(n,\C)$, $Sp(2n,\C)$ and $SO(2n+1,\C)$
with $N=n$, $2n$ and $2n+1$, respectively. Setting $\x=(x_1,x_2,\ldots,x_n)$
and $\ovx=(\ov{x}_1,\ov{x}_2,\ldots,\ov{x}_n)$ with $\ov{x}_i=x_i^{-1}$ for $i=1,2,\ldots,n$, 
for use throughout this paper, we adopt a notation
akin to that used for Schur functions, whereby $s_\lambda(\x)=\ch V_{GL(n,\C)}^\lambda(\x)$,
$sp_\lambda(\x,\ovx)=\ch V_{Sp(2n,\C)}^\lambda(\x,\ovx)$ and $so_\lambda(\x,\ovx,1)=\ch V_{SO(2n+1,\C)}^\lambda(\x,\ovx,1)$.

In the case of $GL(n,\C)$ the transition from a Schur function, $s_\lambda(\x)$, to a factorial Schur function $s_\lambda(\x\,|\,\a)$ 
involves an infinite sequence of factorial parameters $\a=(a_1,a_2,\ldots)$. The transition is effected
by replacing each non-negative power $x_i^m$ of $x_i$ 
by its factorial power $(x_i\,|\,\a)^m$ defined by:
\begin{equation}\label{eqn-pm-fact}
    (x_i\,|\,\a)^m = \begin{cases} (x_i+a_1)(x_i+a_2)\cdots(x_i+a_m) &\hbox{if $m>0$};\cr 1&\hbox{if $m=0$}. \cr\end{cases}
\end{equation}
In the case of the other groups in order to accomodate negative powers $x_i^{-m}=\ov{x}_i^m$ of $x_i$ it is convenient to let
\begin{equation}\label{eqn-pinvm-fact}
    (\ov{x}_i\,|\,\a)^m = \begin{cases} (\ov{x}_i+a_1)(\ov{x}_i+a_2)\cdots(\ov{x}_i+a_m) &\hbox{if $m>0$};\cr 1&\hbox{if $m=0$}. \cr\end{cases}
\end{equation}

We then propose the following definition of factorial characters of the classical Lie groups:
\begin{Definition}\label{Def-fchar}
For any partition $\lambda=(\lambda_1,\lambda_2,\ldots,\lambda_n)$ 
and any $\a=(a_1,a_2,\ldots)$ let
\begin{align}
s_\lambda(\x\,|\,\a) 
&= \frac{\left|\, (x_i\,|\,\a)^{\lambda_j+n-j} \,\right|}
				{\left|\, (x_i\,|\,\a)^{n-j} \,\right| }\,; \label{eqn-gla}\\
sp_\lambda(\x,\ovx\,|\,\a) 
&= \frac{\left|\, x_i(x_i\,|\,\a)^{\lambda_j+n-j} - \ov{x}_i(\ov{x}_i\,|\,\a)^{\lambda_j+n-j} \,\right|}
				                          {\left|\, x_i(x_i\,|\,\a)^{n-j} - \ov{x}_i(\ov{x}_i\,|\,\a)^{n-j} \,\right|}\,; \label{eqn-spa}\\
so_\lambda(\x,\ovx,1\,|\,\a)
&=  \frac{\left|\, x_i^{1/2} (x_i\,|\,\a)^{\lambda_j+n-j} - \ov{x}_i^{1/2}(\ov{x}_i\,|\,\a)^{\lambda_j+n-j} \,\right|}
				                          {\left|\, x_i^{1/2}(x_i\,|\,\a)^{n-j} - \ov{x}_i^{1/2}(\ov{x}_i\,|\,\a)^{n-j} \,\right|}\,. \label{eqn-soa}																
\end{align}
\end{Definition}

Setting $\a=\0=(0,0,\ldots)$ one recovers the classical non-factorial characters of $GL(n,\C)$, $Sp(2n,\C)$ and $SO(2n+1,\C)$,
as given for example in~\cite{FH04}. 
That these definitions are appropriate for non-zero $\a$, in particular the separation out of the factors $x_i$, $\ov{x}_i$, $x_i^{1/2}$
and $\ov{x}_i^{1/2}$,
depends to what extent the properties of these factorial characters are truly analagous to
those of factorial Schur functions. We have in mind things like deriving for each of our factorial characters 
some sort of factorial Jacobi-Trudi identity and a combinatorial interpretation in terms of tableaux as established in the case of factorial Schur 
functions by Macdonald~\cite{Mac92}, and perhaps more ambitiously, the derivation of Tokuyama type identities~\cite{Tok88} 
as recently derived in the factorial Schur function case by Bump, McNamara and Nakasuji~\cite{BMN14}, with an alternative derivation 
appearing in~\cite{HK15}. The key to accomplishing all this is the identification of appropriate analogues of 
the complete homogeneous symmetric functions $h_m(\x)$. This is done for each group in the case of both ordinary and 
factorial characters in Section~\ref{sec-fJT}, in which merely by manipulating determinants, factorial flagged Jacobi-Trudi
identities are derived for each of our factorial characters. 

This is followed in Section~\ref{sec-lambda-m} by consideration of the special cases of one part partitions $\lambda=(m)$.
This enables us to build up in Section~\ref{sec-tabT-LP} a combinatorial realisation of factorial characters, first 
in terms of non-intersecting lattice path models and then in terms of the tableaux traditionally used to specify 
classical, non-factorial characters. Encouraged by this, we present in Section~\ref{sec-tabP-factQ} definitions
of factorial $Q$-functions in terms of appropriate primed shifted tableaux in the general linear, symplectic and odd orthogonal cases.
Reversing our previous trajectory, we proceed in Section~\ref{sec-PtoLP-detQ} from these primed tableaux to non-intersecting lattice path models
and thence to determinantal expressions for each of our factorial $Q$-function. Finally, in Section~\ref{sec-Tok} it is demonstrated
that in the case $\ell(\lambda)=n$ each of our factorial $Q$-functions factorises into a simple product multiplied by
a factorial character, thereby generalising Tokuyama's classical identity~\cite{Tok88} to this factorial context.

\section{Factorial flagged Jacobi-Trudi identities} \label{sec-fJT}

In the study of symmetric functions a key role is played by the complete homogeneous symmetric functions
$h_m(\x)$, with $h_\lambda(\x)=h_{\lambda_1}(\x)h_{\lambda_2}(\x)\ldots h_{\lambda_n}(\x)$ forming
a multiplicative basis of the ring of symmetric functions, in terms of which we have the Jacobi-Trudi identity
$s_\lambda(\x)=|h_{\lambda_{j-j+i}}(\x)|$. 
In order to try to establish factorial Jacobi-Trudi identites we need analogues $h_m(\x\,|\,\a)$ of $h_m(\x)$ that are appropriate not only to 
the case of the other group characters in the case $\a=\0$ 
but also to the case of our factorial characters for general $\a$. 
Just as is done classically for $h_m(\x)$ it is convenient to define all these analogues $h_m(\x\,|\,\a)$ 
by means of generating functions. Our notation is such that each generating function $F(\z,\a;t)$ may be expanded
as a power series in $t$, and we denote the coefficient of $t^m$ in such an expansion by $[t^m]\ F(\z,\a;t)$ for all integers $m$.
In our factorial situation we make the following definitions in terms of generating functions $F_m(\z,\a;t)$ that are truncated 
in the sense that the power $m$ of $[t^m]$ appears in an upper limit of the associated generating function.

\begin{Definition}\label{Def-hma}
For any integer $m\geq0$ and $\a=(a_1,a_2,\ldots)$ let 
\begin{align}
h^{\gl}_m(\x\,|\,\a) &= [t^m]\ \prod_{i=1}^n \frac{1}{1-tx_i} \ \prod_{j=1}^{n+m-1}(1+ta_j)\,;\label{eqn-gl-hma} \\
h^{\sp}_m(\x,\ovx\,|\,\a) &= [t^m]\ \prod_{i=1}^n \frac{1}{(1-tx_i)(1-t\ov{x}_i)}\ \prod_{j=1}^{n+m-1}(1+ta_j)\,;\label{eqn-sp-hma} \\
h^{\oo}_m(\x,\ovx,1\,|\,\a) &= [t^m]\ (1+t)\ \prod_{i=1}^n \frac{1}{(1-tx_i)(1-t\ov{x}_i)}\ \prod_{j=1}^{n+m-1}(1+ta_j)\,;\label{eqn-oo-hma} 
\end{align}
Then for $m=0$ we have $h^{\gl}_0(\x\,|\,\a)= h^{\sp}_0(\x,\ovx\,|\,\a) = h^{\oo}_0(\x,\ovx,1\,|\,\a) = 1$,
while for $m<0$ we set $h^{\gl}_m(\x\,|\,\a)= h^{\sp}_m(\x,\ovx\,|\,\a) = h^{\oo}_m(\x,\ovx,1\,|\,\a) = 0$.
\end{Definition}

The one variable case $\x=(x_i)$ of these Definitions~\ref{Def-hma} allow us to rewrite
our factorial characters in the following manner:
\begin{Lemma}\label{Lem-char-hmxi}
For any partition $\lambda=(\lambda_1,\lambda_2,\ldots,\lambda_n)$ and $\a=(a_1,a_2,\ldots)$ 
\begin{align}
s_\lambda(\x\,|\,\a) 
&= \frac{\left|\, h^{\gl}_{\lambda_j+n-j}(x_i\,|\,\a) \,\right|}
				{\left|\, h^{\gl}_{n-j}(x_i\,|\,\a) \,\right| }\,; \label{eqn-gl-hm}\\						
sp_\lambda(\x,\ovx\,|\,\a) 
&= \frac{\left|\, h^{\sp}_{\lambda_j+n-j}(x_i,\ov{x}_i\,|\,\a) \,\right|}
				                          {\left|\, h^{\sp}_{n-j}(x_i,\ov{x_i}\,|\,\a) \,\right|}\,; \label{eqn-sp-hm}\\
so_\lambda(\x,\ovx,1\,|\,\a)
&=  \frac{\left|\, h^{\oo}_{\lambda_j+n-j}(x_i,\ov{x}_i,1\,|\,\a) \,\right|}
				                          {\left|\, h^{\oo}_{n-j}(x_i,\ov{x}_i\,|\,\a) \,\right|}\,. \label{eqn-so-hm}
\end{align}
\end{Lemma}

\noindent{\bf Proof}:
In the case of 
$s_\lambda(\x\,|\,\a)$ it suffices to note that for $m\geq 0$
\begin{align}\label{eqn-gl-hxma}
    h_m(x_i\,|\,\a) &= [t^m] \frac{1}{1-tx_i} \prod_{j=1}^m (1+ta_j) =  [t^m] \frac{1+ta_m}{1-tx_i} \prod_{j=1}^{m-1} (1+ta_j)\cr
                &= [t^m] \left(1+\frac{t(x_i+a_m)}{1-tx_i}\right) \prod_{j=1}^{m-1} (1+ta_j) = (x_i+a_m) [t^{m-1}] \frac{1}{1-tx_i} \prod_{j=1}^{m-1} (1+ta_j)\cr
								&= (x_i+a_m)(x_i+a_{m-1})\cdots(x_i+a_1) [t^0] \frac{1}{1-tx_i} = (x_i\,|\,\a)^m\,.
\end{align}
One then just uses this identity in (\ref{eqn-gla}) with $m=\lambda_j+n-j$ and $m=n-j$
in the numerator and denominator, respectively.

In the case $sp_\lambda(\x,\ovx\,|\,\a)$ we have
\begin{align}\label{eqn-sp-hxma}
h_m^{\sp}(x_i,\ov{x}_i\,|\,\a)&=[t^m]\frac{1}{(1-tx_i)(1-t\ov{x}_i)} \prod_{j=1}^m (1+ta_j) \cr
	          &=[t^m] \frac{1}{x_i-\ov{x}_i}\left(\frac{x_i}{1-tx_i}-\frac{\ov{x}_i}{1-t\ov{x}_i}\right) \prod_{j=1}^m (1+ta_j)\cr
	                    &=\frac{1}{x_i-\ov{x}_i}\left( x_i(x_i\,|\,\a)^m-\ov{x}_i(\ov{x}_i\,|\,\a)^m\right)\,,
\end{align}
where use has been made of (\ref{eqn-gl-hxma}).
The required result follows by using these identities in (\ref{eqn-spa}) as before,
with the cancellation between numerator and denominator of the common factors $x_i-\ov{x}_i$ for $i=1,2,\ldots,n$.

Similarly, in the $so_\lambda(\x,\ovx,1\,|\,\a)$ case we have
\begin{align}\label{eqn-so-hxma}
h_m^{\oo}(x_i,\ov{x}_i,1\,|\,\a)&=[t^m]\frac{1+t}{(1-tx_i)(1-t\ov{x}_i)} \prod_{j=1}^m (1+ta_j)\cr
	          &=[t^m] \frac{1}{x_i^{1/2}-\ov{x}_i^{1/2}}\left(\frac{x_i^{1/2}}{1-tx_i}-\frac{\ov{x}_i^{1/2}}{1-t\ov{x}_i}\right) \prod_{j=1}^m (1+ta_j)\cr
	                    &=\frac{1}{x_i^{1/2}-\ov{x}_i^{1/2}}\left( x_i^{1/2}(x_i\,|\,\a)^m-\ov{x}_i^{1/2}(\ov{x}_i\,|\,\a)^m\right)\,.
\end{align}
Then one again uses this identity in (\ref{eqn-soa}) as before,
with the cancellation this time of the common factors $x_i^{1/2}-\ov{x}_i^{1/2}$ for $i=1,2,\ldots,n$.
\qed

The next step is to transform each of the expressions in Lemma~\ref{Lem-char-hmxi} into some sort of flagged Jacobi-Trudi identity.
This is accomplished by means of the following Lemma:

\begin{Lemma}\label{Lem-hra-xixj}
For all $i$ and $j$ such that $1\leq i<j\leq n$ and all integers $m$:
\begin{align}
  &h^{\gl}_m(x_i,\ldots x_{j-1}\,|\,\a)-h^{\gl}_m(x_{i+1},\ldots,x_j\,|\,\a)\cr
	&~~~~~~~~~~~~=(x_i-x_j) h^{\gl}_{m-1}(x_i,\ldots,x_j\,|\,\a)\,; \label{eqn-gl-hra-xixj} \\
	&h^{\sp}_m(x_i,\ov{x}_i,\ldots,x_{j-1},\ov{x}_{j-1}\,|\,\a)-h^{\sp}_m(x_{i+1},\ov{x}_{i+1},\ldots,x_{j},\ov{x}_{j},\,|\,\a)\cr
	&~~~~~~~~~~~~=(x_i-x_j)(1-\ov{x}_i\ov{x}_j) h^{\sp}_{m-1}(x_i,\ov{x}_i,\ldots,x_{j},\ov{x}_{j},\,|\,\a);\; \label{eqn-sp-hra-xixj} \\
	&h^{\oo}_m(x_i,\ov{x}_i,\ldots,x_{j-1},\ov{x}_{j-1},1\,|\,\a)-h^{\oo}_m(x_{i+1},\ov{x}_{i+1},\ldots,x_{j},\ov{x}_{j},1\,|\,\a)\cr
	&~~~~~~~~~~~~=(x_i-x_j)(1-\ov{x}_i\ov{x}_j) h^{\oo}_{m-1}(x_i,\ov{x}_i,\ldots,x_{j},\ov{x}_{j},1\,|\,\a)\,. \label{eqn-so-hra-xixj} 
\end{align}
\end{Lemma}

\noindent{\bf Proof}: First it should be noted that all these identities are trivially true for $m<0$ and for $m=0$ since
each $h_{m}$ reduces to either $0$ or $1$, with each $h_{m-1}$ reducing to $0$.
For $m>0$, in the simplest case
\begin{align}
     &h^{\gl}_m(x_i,\ldots x_{j-1}\,|\,\a)\!-\!h^{\gl}_m(x_{i+1},\ldots,x_j\,|\,\a)
		 =[t^m]\ \left((1\!-\!tx_j)\!-\!(1\!-\!tx_i)\right)\prod_{\ell=i}^j \frac{1}{1\!-\!tx_\ell}\!\!\!\prod_{k=1}^{m+j-i-1}\!\!(1\!+\!ta_k)\cr
	   &= (x_i\!-\!x_j)\ [t^{m-1}]\  \prod_{\ell=i}^j \frac{1}{1\!-\!tx_\ell} \!\!\!\prod_{k=1}^{(m-1)+j-i}\!\!\! (1\!+\!ta_k) 
		 = (x_i\!-\!x_j)h^{\gl}_{m-1}(x_i,\ldots,x_j\,|\,\a)\,. \label{eqn-gl-hrHa-xixj}
\end{align}

The other two cases are essentially the same, and are illustrated by the symplectic case:
\begin{align}
  &h^{\sp}_m(x_i,\ov{x}_i,\ldots,x_{j-1},\ov{x}_{j-1}\,|\,\a)-h^{\sp}_m(x_{i+1},\ov{x}_{i+1},\ldots,x_{j},\ov{x}_{j},\,|\,\a)\cr
&= [t^m]\ \left( (1\!-\!tx_j)(1\!-\!t\ov{x}_j)\!-\!(1\!-\!tx_i)(1\!-\!t\ov{x}_i)\right)\ \prod_{\ell=i}^j \frac{1}{(1\!-\!tx_\ell)(1\!-\!t\ov{x}_\ell)} 
          \!\!\prod_{k=1}^{m+j-i-1}\!\!(1+ta_k)\cr
		 &= (x_i\!+\!\ov{x}_i\!-\!x_j\!-\!\ov{x}_j)\ [t^{m-1}]\ \prod_{\ell=i}^j \frac{1}{(1\!-\!tx_\ell)(1\!-\!t\ov{x}_\ell)} 
		       \!\!\prod_{k=1}^{(m-1)+j-i}\!\!(1+ta_k)\cr 
	   &= (x_i\!-\!x_j)(1\!-\!\ov{x}_i\ov{x}_j)\ h^{\sp}_{m-1}(x_i,\ov{x}_i,\ldots,x_{j},\ov{x}_{j},\,|\,\a)\,. \label{eqn-sp-hrHa-xixj}
\end{align}
Exactly the same procedure applies to the odd orthogonal case.
\qed

Now we are a position to state and prove the following result:

\begin{Theorem}\label{The-fJT} 
Let $\x=(x_1,x_2,\ldots,x_n)$ and $\ovx=(\ov{x}_1,\ov{x}_2,\ldots,\ov{x}_n)$ with
$\x^{(i)}=(x_i,x_{i+1},\ldots,x_n)$ and $\ovx^{(i)}=(\ov{x}_i,\ov{x}_{i+1},\ldots,\ov{x}_n)$ 
for $i=1,2,\ldots,n$. 
Then for any partition $\lambda$ of length $\ell(\lambda)\leq n$ and any $\a=(a_1,a_2,\ldots)$ we have 
\begin{align}
s_\lambda(\x\,|\,\a) 
&= \left|\, h^{\gl}_{\lambda_j-j+i}(\x^{(i)}\,|\,\a) \,\right|; \label{eqn-gl-fJT} \\						
sp_\lambda(\x,\ovx\,|\,\a) 
&= \left|\, h^{\sp}_{\lambda_j-j+i}(\x^{(i)},\ovx^{(i)}\,|\,\a) \,\right|; \label{eqn-sp-fJT}\\
so_\lambda(\x,\ovx,1\,|\,\a)
&=  \left|\, h^{\oo}_{\lambda_j-j+i}(\x^{(i)},\ovx^{(i)},1\,|\,\a) \,\right|. \label{eqn-so-fJT}
\end{align}
\end{Theorem}

\noindent{\bf Proof}: In the classical non-factorial case, obtained by setting $\a=\0$, these flagged Jacobi-Trudi identities have all been obtained previously
by Okada~\cite{Oka??} by means of lattice path methods. The symplectic case was also independently obtained by this means
by Chen, Li and Louck~\cite{CLL02}, while the Schur function case goes back at least to Littlewood~\cite{Lit50} 
who obtained it en route to his derivation of the classical Jacobi-Trudi identity by means of simple row manipulations of determinants.

It is Littlewood's method that we use here to establish all three identities. 
Subtracting row $(i+1)$ from row $i$ for $i=1,2,\ldots,n-1$ in the numerator of (\ref{eqn-gl-hm}) and applying (\ref{eqn-gl-hra-xixj}), 
then repeating the process for $i=1,2,\ldots,n-2$ and so on, yields
\begin{align}\label{eqn-gl-hlambda}
& \left|\, h_{\lambda_j+n-j}(x_i\,|\,\a) \,\right| 
= \left|\,\begin{array}{c} h_{\lambda_j+n-j}(x_i\,|\,\a) - h_{\lambda_j+n-j}(x_{i+1}\,|\,\a) \cr
                                                                                  h_{\lambda_j+n-j}(x_n\,|\,\a) \cr 
																																									 \end{array}    \,\right| \cr
& = \prod_{i=1}^{n-1} (x_i-x_{i+1}) \ \left|\,\begin{array}{c} h_{\lambda_j+n-j-1}(x_i,x_{i+1}\,|\,\a) \cr
                                                                            h_{\lambda_j+n-j}(x_n\,|\,\a) \cr 
																																									 \end{array}  \,\right| \cr  
&= \prod_{i=1}^{n-1} (x_i-x_{i+1}) \ \left|\,\begin{array}{c} h_{\lambda_j+n-j-1}(x_i,x_{i+1}\,|\,\a)-h_{\lambda_j+n-j-1}(x_{i+1},x_{i+2} \,|\,\a) \cr
                                                                          h_{\lambda_j+n-j-1}(x_{n-1},x_n\,|\,\a) \cr 
																																						h_{\lambda_j+n-j}(x_n\,|\,\a) \cr 
    																																							 \end{array}  \,\right| \cr																																
& = \prod_{i=1}^{n-1} (x_i-x_{i+1})\prod_{i=1}^{n-2} (x_i-x_{i+2}) \ \left|\,\begin{array}{c} h_{\lambda_j+n-j-2}(x_i,x_{i+1},x_{i+2}\,|\,\a) \cr
                                                                           h_{\lambda_j+n-j-1}(x_{n-1},x_n\,|\,\a) \cr 
                                                                            h_{\lambda_j+n-j}(x_n\,|\,\a) \cr 
																																									 \end{array}  \,\right| \cr  
& = ~~\cdots~~ =   \prod_{1\leq i<j\leq n} (x_i-x_j)\ \left|\, h_{\lambda_j+n-j-(n-i)}(x_i,x_{i+1},\ldots,x_n\,|\,\a) \,\right| \cr
& = \prod_{1\leq i<j\leq n} (x_i-x_j)\ \left|\, h_{\lambda_j-j+i}(\x^{(i)}\,|\,\a) \,\right|\,. 
\end{align}
In the special case $\lambda=(0)$ this yields the denominator identity
\begin{equation}\label{eqn-gl-h0}
         \left|\, h_{n-j}(x_i\,|\,\a) \,\right|=\prod_{1\leq i<j\leq n} (x_i-x_j)\ \left|\, h_{-j+i}(\x^{(i)}\,|\,\a) \,\right|=\prod_{1\leq i<j\leq n} (x_i-x_j)
\end{equation} 
since the determinant $\left|\, h_{-j+i}(\x^{(i)}\,|\,\a) \,\right|$ is lower-triangular with all its diagonal elements equal to $h_0(\x^{(i)}\,|\,\a)=1$.
Taking the ratio of these two formulae implies that $s_\lambda(x\,|\,\a) =  \left|\, h_{\lambda_j-j+i}(\x^{(i)}\,|\,\a) \,\right|$, 
as required in (\ref{eqn-gl-fJT}).

The other two cases (\ref{eqn-sp-fJT}) and (\ref{eqn-so-fJT}) are obtained in an identical manner. 
The only difference is that instead of extracting factors $(x_i-x_j)$ as dictated by (\ref{eqn-gl-hra-xixj}), 
one extracts factors $(x_i-x_j)(1-\ov{x}_i\ov{x}_j)$ as dictated by both (\ref{eqn-sp-hra-xixj}) and (\ref{eqn-so-hra-xixj}). 
\qed

\section{Explicit formulae in the case $\lambda=(m)$}\label{sec-lambda-m}

As a consequence of Theorem~\ref{The-fJT} it should be noted that in the 
case of a one-part partition $\lambda=(m,0,\ldots,0)=(m)$ we have
\begin{Corollary}
For all non-negative integers $m$
\begin{equation}
s_{(m)}(\x\,|\,\a) = h^{\gl}_{m}(\x\,|\,\a);~~ 
sp_{(m)}(\x,\ovx\,|\,\a) = h^{\sp}_{m}(\x,\ovx\,|\,\a);~~
so_{(m)}(\x,\ovx,1\,|\,\a) = h^{\oo}_{m}(\x,\ovx,1\,|\,\a) \,. \label{eqn-corr} 
\end{equation}
\end{Corollary}

\noindent{\bf Proof}:~~On setting $\lambda=(m,0,\ldots,0)$ the flagged factorial Jacobi-Trudi determinants in (\ref{eqn-gl-fJT})-(\ref{eqn-so-fJT}) 
are reduced to lower-triangular form since each $h_{-j+i}=0$ for $i<j$.  Moreover for $i>1$ the diagonal entries are
all of the form $h_0=1$, while the $(1,1)$ entry is just $h_{m}$ with $\x^{(1)}=\x$ and $\ovx^{(1)}=\ovx$. 
\qed

Factorial characters in the one-part partition case may then be evaluated directly from the generating function 
formulae of Definition~\ref{Def-hma}. Before doing so it is convenient, following Macdonald~\cite{Mac92},
to introduce the shift operator $\tau$ defined by
\begin{equation}\label{eqn-tau}
    \tau^r \a = (a_{r+1},a_{r+2},\ldots) \quad\hbox{for any integer $r$ and any $\a=(a_1,a_2,\ldots)$}\,.
\end{equation}	

In the Schur function case with $\x=(x_1,x_2,\ldots,x_n)$ and $\x'=(x_1,x_2,\ldots,x_{n-1})$ the 
generating function (\ref{eqn-gl-hm}) yields
\begin{align}\label{eqn-hm-recur}
  h_m(\x\,|\,\a) &= [t^m]\ \prod_{i=1}^n \frac{1}{1-tx_i}\  \prod_{k=1}^{m+n-1}(1+ta_k)\cr
	               &= [t^m]\ \left( \frac{1+ta_{m+n-1}}{1-tx_n}\right)\ \prod_{i=1}^{n-1} \frac{1}{1-tx_i}\  \prod_{k=1}^{m+n-2}(1+ta_k)\cr
								 &= [t^m]\ \left( 1+ \frac{t(x_n+a_{n+m-1})}{1-tx_n}\right)\ \prod_{i=1}^{n-1} \frac{1}{1-tx_i}\  \prod_{k=1}^{m+n-2}(1+ta_k)\cr
								 &= h_m(\x'\,|\,\a)+(x_n+a_{m+n-1}) h_{m-1}(\x\,|\,\a)\,.
\end{align}
Iterating this recursion relation gives
\begin{equation}\label{eqn-gl-hm-expansion}
  h_m(\x\,|\,\a)= \sum_{1\leq i_1\leq i_2\leq \cdots\leq i_m\leq n} (x_{i_1}+a_{i_1})(x_{i_2}+a_{i_2+1})\cdots(x_{i_m}+a_{i_m+m-1})\,.
\end{equation}

This result can be exploited in the symplectic case, where it might be noted first that if we introduce dummy
parameters $a_\ell=0$ for $\ell=0,-1,-2,\ldots$ then it follows from (\ref{eqn-sp-hm}) that
\begin{align}
 & h^{\sp}_m(\x,\ovx\,|\,\a) = [t^m]\ \prod_{i=1}^n \frac{1}{(1-tx_i)(1-t\ov{x}_i)}\  \prod_{k=1}^{m+n-1}(1+ta_k)\cr
 &= [t^m]\ \prod_{i=1}^n \frac{1}{(1-tx_i)(1-t\ov{x}_i)}\!\!\prod_{k=1-n}^{m+2n-1-n}\!\!(1+ta_k)\!
													 = h_m(\x,\ovx\,|\,\tau^{-n}\a) = h_m(\z\,|\,\tau^{-n}\a) \label{eqn-sp-hmtau}
\end{align}
where $\tau^{-n}\a=(a_{-n+1},\ldots,a_{-1},a_0,a_1,a_2\ldots)$, and it is convenient to order the indeterminates in $\z$
so that $\z=(x_1,\ov{x}_1,x_2,\ov{x}_2,\ldots,x_n,\ov{x}_n)$. It then follows that
\begin{equation}\label{eqn-sp-hm-expansion}
  h^{\sp}_m(\x,\ovx\,|\,\a)= \sum_{1\leq i_1\leq i_2\leq \cdots\leq i_m\leq 2n} (z_{i_1}+a_{i_1-n})(z_{i_2}+a_{i_2-n+1})\cdots(z_{i_m}+a_{i_m-n+m-1})\,.
\end{equation}
with 
\begin{equation}
       z_{i_j}+a_{i_j-n+j-1}=\begin{cases} x_k+a_{2k-n+j-2}&\hbox{if $i_j=2k-1$;}\cr
			                                     \ov{x}_k+a_{2k-n+j-1}&\hbox{if $i_j=2k$,}\cr
																					\end{cases} ~~~~\hbox{with $a_\ell=0$ if $\ell\leq 0$}.
\end{equation}

Turning to the odd orthogonal case and using (\ref{eqn-so-hm}) we have
\begin{align}
  &h^{\oo}_m(\x,\ovx,1\,|\,\a) = [t^m]\ (1+t)\ \prod_{i=1}^n \frac{1}{(1-tx_i)(1-t\ov{x}_i)}\  \prod_{k=1}^{m+n-1}(1+ta_k)\cr
	                         &= [t^m]\ \left((1+ta_{m+n})+t(1-a_{m+n})\right) \prod_{i=1}^n \frac{1}{(1-tx_i)(1-t\ov{x}_i)}\  \prod_{k=1}^{m+n}(1+ta_k)\cr
													&= h^{\sp}_m(\x,\ovx\,|\,\tau\a) + (1-a_{m+n})\ h^{\sp}_{m-1}(\x,\ovx\,|\,\tau\a) \cr 
													&= h_m(\x,\ovx\,|\,\tau^{1-n} \a) + (1-a_{m+n})\ h_{m-1}(\x,\ovx\,|\,\tau^{1-n}\a)\,, \label{eqn-so-hmtau}
\end{align}
where dummy parameters $a_\ell=0$ for $\ell=0,-1,-2,\ldots$ are once again involved.

It follows that
\begin{align}\label{eqn-so-hm-expansion}
  &h^{\oo}_m(\x,\ovx,1\,|\,\a)= \sum_{1\leq i_1\leq i_2\leq \cdots\leq i_m\leq 2n}\hskip-2ex(z_{i_1}\!+\!a_{i_1+1-n})(z_{i_2}\!+\!a_{i_2+2-n})\cdots(z_{i_m}\!+\!a_{i_m-n+m})\cr
  &+\hskip-1ex\sum_{1\leq i_1\leq i_2\leq \cdots\leq i_{m-1}\leq 2n}\hskip-4ex(z_{i_1}\!+\!a_{i_1+1-n})(z_{i_2}\!+\!a_{i_2+2-n})\cdots(z_{i_{m-1}}\!+\!a_{i_{m-1}+m-1-n})(1\!-\!a_{m+n})
\end{align}
with 
\begin{equation}
       z_{i_j}+a_{i_j-n+j}=\begin{cases} x_k+a_{2k-n+j-1}&\hbox{if $i_j=2k-1$;}\cr
			                                     \ov{x}_k+a_{2k-n+j}&\hbox{if $i_j=2k$,}\cr
																					\end{cases} ~~~~\hbox{with $a_\ell=0$ if $\ell\leq 0$}.
\end{equation}

\section{Combinatorial realisation of factorial characters}\label{sec-tabT-LP}

The significance of these results is that they offer an immediate lattice path model of each of the relevant
one-part partition factorial characters. Then by making use of $n$-tuples of such lattice paths in the interpretation
of the factorial flagged Jacobi-Trudi identities of Theorem~\ref{The-fJT} one arrives at a non-intersecting
lattice path model of factorial characters specified by any partition $\lambda$ of length $\ell(\lambda)\leq n$. 
This leads inexorably to a further realisation of factorial characters in terms of certain
appropriately weighted tableaux. The tableaux themselves are none other than those already associated with
Schur functions, symplectic group characters and odd orthogonal group characters in the classical non-factorial 
case.

Restricting our attention to fixed $n$ and partitions
$\lambda=(\lambda_1,\lambda_2,\ldots,\lambda_n)$ of length $\ell(\lambda)\leq n$,
each such partition defines a Young diagram $F^\lambda$
consisting of $|\lambda|=\lambda_1+\lambda_2+\cdots+\lambda_n$ boxes arranged in $\ell(\lambda)$ rows of lengths $\lambda_i$
for $i=1,2,\ldots,\ell(\lambda)$. We adopt the English convention as used by Macdonald~\cite{Mac95} whereby the rows are
left-adjusted to a vertical line and are weakly decreasing in length from top to bottom.
For example in the case $\lambda=(4,3,3)$, for which $\ell(\lambda)=3$ and
$$
F^{433}\ =\ 
\YT{0.2in}{0.2in}{}{
 {{},{},{},{}},
 {{},{},{}},
 {{},{},{}},
}
$$
More precisely we define $F^\lambda=\{ (i,j)\,|\, 1\leq i\leq \ell(\lambda) ; 1\leq j\leq \lambda_i\}$ and refer
to $(i,j)$ as being the box in the $i$th row and $j$th column of $F^\lambda$. Assigning an entry $T_{ij}$
taken from some alphabet to each box $(i,j)$ of $F^\lambda$ in accordance with various rules gives rise to 
tableaux $T$ of shape $\lambda$ that may be used, as we shall see, to express both ordinary and factorial characters 
in a combinatorial manner.

\begin{Definition} [Littlewood~\cite{Lit37}]~\label{Def-ssT}
Let ${\cal T}^{\gl}_\lambda$ be the set of all semistandard Young tableaux $T$ of shape $\lambda$ 
that are obtained by filling each box $(i,j)$ of $F^\lambda$ with an entry $T_{ij}$ from the alphabet
$$
     \{ 1<2<\cdots<n \}
$$
in such a way that: 
{\bf (T1)} entries weakly increase across rows from left to right;
{\bf (T2)} entries weakly increase down columns from top to bottom;
{\bf (T3)} no two identical non-zero entries appear in the same column.
\end{Definition}

\begin{Definition} [King~\cite{Kin76}]~\label{Def-spT}
Let ${\cal T}^{\sp}_\lambda$ be the set of all symplectic tableaux $T$ of shape $\lambda$
that are obtained by filling each box $(i,j)$ of $F^\lambda$ with an entry $T_{ij}$ from
the alphabet
$$
    \{1<\ov{1}<2<\ov{2}<\cdots<n<\ov{n}\}
$$
in such a way that conditions {\bf (T1)-(T3)} are satisfied, together with: 
{\bf (T4)} neither $k$ nor $\ov{k}$ appear lower than the $k$th row.
\end{Definition}

\begin{Definition} [Sundaram~\cite{Sun90}]~\label{Def-soT}
Let ${\cal T}^{\oo}_\lambda$ be the set of all odd orthogonal tableaux $T$ of shape $\lambda$
that are obtained by filling each box $(i,j)$ of $F^\lambda$ with an entry $T_{ij}$ from
the alphabet
$$
    \{1<\ov{1}<2<\ov{2}<\cdots<n<\ov{n}<0\}
$$
in such a way that conditions {\bf (T1)-(T4)} are satisfied, together with: 
{\bf (T5)} in any row $0$ appears at most once.
\end{Definition}

These definitions are exemplified for $GL(4)$, $Sp(8)$ and $SO(9)$ 
in turn as shown below from left to right:
\begin{equation}\label{eqn-tabx3}
\YT{0.2in}{0.2in}{}{
 {{1},{1},{2},{4}},
 {{2},{3},{3}},
 {{4},{4},{4}},
} 
\qquad
\YT{0.2in}{0.2in}{}{
 {{1},{\ov1},{2},{\ov4}},
 {{\ov3},{4},{4}},
 {{4},{\ov4},{\ov4}},
}
\qquad
\YT{0.2in}{0.2in}{}{
 {{1},{\ov1},{2},{\ov4}},
 {{\ov3},{4},{0}},
 {{4},{\ov4},{0}},
}
\end{equation}

These definitions allow us to provide combinatorial expressions for factorial characters as follows:
\begin{Theorem}\label{The-factchar-tab}
For each $\g$ and $\z$ as tabulated below, and any $\a=(a_1,a_2,\ldots)$
\begin{equation}\label{eqn-factchar-tab}
   \g_\lambda(\z\,|\,\a) = \sum_{T\in{\cal T}^{\g}_\lambda}\ \prod_{(i,j)\in F^\lambda} \wgt(T_{ij})\,.
\end{equation}
where
\begin{equation}\label{eqn-Twgts}
\begin{array}{|l|l|l|l|}
\hline
\g_\lambda(\z\,|\,\a)&T_{ij}&\wgt(T_{ij})&\cr
\hline
s_\lambda(\x\,|\,\a)&k &x_k+a_{k+j-i}&\cr
\hline
sp_\lambda(\x,\ovx\,|\,\a)&k&x_k+a_{2k-1-n+j-i}&a_m=0\hbox{~~for~~}m\leq0\cr
                     &\ov{k}&\ov{x}_k+a_{2k-n+j-i}&\cr
\hline
so_\lambda(\x,\ovx,1\,|\,\a)&k&x_k+a_{2k-n+j-i}&a_m=0\hbox{~~for~~}m\leq0\cr
                     &\ov{k}&\ov{x}_k+a_{2k+1-n+j-i}&\cr
										 &0&1-a_{n+1+j-i}&\cr
\hline	
\end{array}
\end{equation}
with $\x=(x_1,x_2,\ldots,x_n)$, $\ovx=(\ov{x}_1,\ov{x}_2,\ldots,\ov{x}_n)$ and $\lambda$
any partition of length $\ell(\lambda)\leq n$.
\end{Theorem}

\begin{figure*}[htbp]
$$
\begin{array}{c}
LP(T)\ =\ 
\vcenter{\hbox{
\begin{tikzpicture}[x={(0in,-0.3in)},y={(0.3in,0in)}] 
\foreach \j in {4,...,8} \draw(1,\j)node{$\bullet$};
\foreach \j in {3,...,8} \draw(2,\j)node{$\bullet$};
\foreach \j in {2,...,8} \draw(3,\j)node{$\bullet$};
\foreach \j in {1,...,8} \draw(4,\j)node{$\bullet$};
\draw[-](0+0.2,6-0.2)to(4,2);
\draw[-](0+0.2,7-0.2)to(4,3);
\draw[-](0+0.2,8-0.2)to(4,4);
\draw[-](0+0.2,9-0.2)to(4,5);
\draw[-](1+0.2,9-0.2)to(4,6);
\draw[-](2+0.2,9-0.2)to(4,7);
\draw[-](3+0.2,9-0.2)to(4,8);
\draw(1-0.2,3.5)node{$P_1$};
\draw(2-0.2,2.5)node{$P_2$};
\draw(3-0.2,1.5)node{$P_3$};
\draw(4-0.2,0.5)node{$P_4$};
\draw(4.5,1)node{$Q_{4}$};  
\draw(4.5,5)node{$Q_{3}$};
\draw(4.5,6)node{$Q_{2}$};
\draw(4.5,8)node{$Q_{1}$};
\draw(0,6)node{$a_1$};
\draw(0,7)node{$a_2$};
\draw(0,8)node{$a_3$};
\draw(0,9)node{$a_4$};
\draw(1,9)node{$a_5$};
\draw(2,9)node{$a_6$};
\draw(3,9)node{$a_7$};
\draw[draw=magenta,ultra thick] (1,4)to(1,5); \draw(1-0.3,5-0.3)node{$\magenta{x_1}$};
\draw[draw=magenta,ultra thick] (1,5)to(1,6); \draw(1-0.3,6-0.3)node{$\magenta{x_1}$};
\draw[draw=magenta,ultra thick] (1,6)to(2,6); 
\draw[draw=magenta,ultra thick] (2,6)to(2,7); \draw(2-0.3,7-0.3)node{$\magenta{x_2}$};
\draw[draw=magenta,ultra thick] (2,7)to(4,7);
\draw[draw=magenta,ultra thick] (4,7)to(4,8); \draw(4-0.3,8-0.3)node{$\magenta{x_4}$};
\draw[draw=blue,ultra thick] (2,3)to(2,4); 
\draw[draw=blue,ultra thick] (2,4)to(3,4); \draw(2-0.3,4-0.3)node{$\blue{x_2}$};
\draw[draw=blue,ultra thick] (3,4)to(3,5); \draw(3-0.3,5-0.3)node{$\blue{x_3}$};
\draw[draw=blue,ultra thick] (3,5)to(3,6); \draw(3-0.3,6-0.3)node{$\blue{x_3}$};
\draw[draw=blue,ultra thick] (3,6)to(4,6);
\draw[draw=red,ultra thick] (3,2)to(4,2); 
\draw[draw=red,ultra thick] (4,2)to(4,3); \draw(4-0.3,3-0.3)node{$\red{x_4}$};
\draw[draw=red,ultra thick] (4,3)to(4,4); \draw(4-0.3,4-0.3)node{$\red{x_4}$};
\draw[draw=red,ultra thick] (4,4)to(4,5); \draw(4-0.3,5-0.3)node{$\red{x_4}$};
\end{tikzpicture}
}}\cr\cr
T\ = \ \YT{0.2in}{0.2in}{}{
 {\magenta{1},\magenta{1},\magenta{2},\magenta{4}},
 {\blue{2},\blue{3},\blue{3}},
 {\red{4},\red{4},\red{4}},
}
\quad
\wgt(T)\ =\ 
\wideYT{0.2in}{0.6in}{}{
 {\magenta{x_1\!+\!a_1},\magenta{x_1\!+\!a_2},\magenta{x_2\!+\!a_4},\magenta{x_4\!+\!a_7}},
 {\blue{x_2\!+\!a_1},\blue{x_3\!+\!a_3},\blue{x_3\!+\!a_4}},
 {\red{x_4\!+\!a_2},\red{x_4\!+\!a_3},\red{x_4\!+\!a_4}},
} \nonumber
\end{array}
$$
\caption{Contribution to $s_{433}(\x\,|\,\a)$ from $T$ and $LP(T)$.}
\label{fig-gl-TtoLP}
\end{figure*}
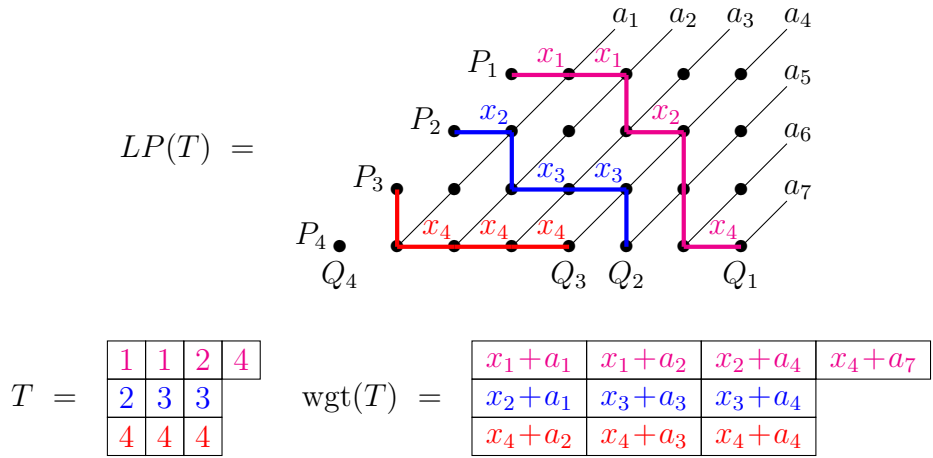

\noindent{\bf Proof}:
In the Schur function case, as in \cite{HK15}, we adopt matrix coordinates $(k,\ell)$ for lattice points with $k=1,2,\ldots,n$ specifying row labels from
top to bottom, and $\ell=1,2,\ldots,\lambda_1+n$ specifying column labels from left to right. 
Each lattice path that we are interested in is a continuous path from some 
$P_i=(i,n-i+1)$ to some $Q_j=(n,n-j+1+\lambda_j)$ with $i,j\in\{1,2,\ldots,n\}$. 
Such a path consists of a sequence of horizontal or vertical edges
and is associated with a contribution to $h_{\lambda_j-j+i}(\x^{(i)}\,|\,\a)$ in the form of
a summand of (\ref{eqn-gl-hma}) with $m=\lambda_j-j+i$ and $\x$ replaced by $\x^{(i)}$. 
Taking into account the restriction of the alphabet from $\x$ to $\x^{(i)}$, the weight assigned to  
horizontal edge from $(k,\ell-1)$ to $(k,\ell)$ is $x_k+a_{k+\ell-n-1}$. 
Thanks to the Lindstr\"om-Gessel-Viennot theorem~\cite{Lin73,GV85,GV89} the only
surviving contributions to the determinantal expression for $s_\lambda(\x\,|\,\a)$ in the flagged factorial Jacobi-Trudi 
identity (\ref{eqn-gl-fJT}) are those corresponding to an $n$-tuple of non-intersecting lattice paths from $P_i$ to $Q_i$ for $i=1,2,\ldots,n$. 
Such $n$-tuples are easily seen~\cite{SW86} to be in bijective correspondence
with semistandard Young tableaux $T$ of shape $\lambda$ as in Definition~\ref{Def-ssT}, with the $j$th horizontal edge at level $k$
on the path from $P_i$ to $Q_i$ giving an entry $T_{ij}=k$ in $T$ for $i=1,2,\ldots,n$ and $j=1,2,\ldots,\lambda_i$. To complete the proof 
of Theorem~\ref{The-factchar-tab} in the factorial Schur function case 
it only remains to note that the weight $\wgt(T_{ij})$ to be assigned to $T_{ij}$ is that of the edge from $(k,\ell-1)$ to $(k,\ell)$  
given by $x_k+a_{k+\ell-n-1}=x_k+a_{k+j-i}$ with $j=\ell-(n-i+1)$ since this is the number of horizontal steps from $P_i$ to column $\ell$
on the lattice path from $P_i$ to $Q_i$.
This is exemplified in Figure~\ref{fig-gl-TtoLP} in the case $n=4$ and $\lambda=(4,3,3)$.

Thanks to (\ref{eqn-sp-hmtau}), the lattice path proof in the factorial symplectic case 
proceeds exactly as in the Schur function case
with the alphabet extended to include both $x_k$ and $\ov{x}_k$ for $k=1,2,\ldots,n$, and with $\a$ replaced by $\tau^{-n}\a$. 
The starting points are now $P_i=(2i-1,n-i+1)$ thereby ensuring that condition {\bf(T4)} is satisfied, 
and the end points are $Q_j=(2n,n-j+1+\lambda_j)$. Once again
it is only the $n$-tuples of non-intersecting lattice paths from $P_i$ to $Q_i$ that
contribute to $sp_\lambda(\x,\ovx\,|\,\a)$ and these are in bijective correspondence with
the symplectic tableaux of Definition~\ref{Def-spT} of shape $\lambda$ with entries from $\{1<\ov1<\cdots<n<\ov{n}\}$.
This is exemplified in Figure~\ref{fig-sp-TtoLP} for $n=4$ and $\lambda=(4,3,3)$. 

\begin{figure*}[htbp]
$$
\begin{array}{c}
LP(T)\ = \ 
\vcenter{\hbox{
\begin{tikzpicture}[x={(0in,-0.3in)},y={(0.3in,0in)}] 
\draw[-](0+0.2,5-0.2)to(1,4);
\draw[-](0+0.2,6-0.2)to(3,3);
\draw[-](0+0.2,7-0.2)to(5,2);
\draw[-](0+0.2,8-0.2)to(7,1);
\foreach \i in {0,...,7} \draw[-,thin](\i+0.2,9-0.2)to(8,\i+1);
\draw(1-0.2,4-0.5)node{$P_1$};
\draw(3-0.2,3-0.5)node{$P_2$};
\draw(5-0.2,2-0.5)node{$P_3$};
\draw(7-0.2,1-0.5)node{$P_4$};
\draw(8+0.5,1)node{$Q_4$};
\draw(8+0.5,5)node{$Q_3$};
\draw(8+0.5,6)node{$Q_2$};
\draw(8+0.5,8)node{$Q_1$};
\draw(0,5)node{$a_{\ov4}$};
\draw(0,6)node{$a_{\ov3}$};
\draw(0,7)node{$a_{\ov2}$};
\draw(0,8)node{$a_{\ov1}$};
\foreach \i in {0,...,7} \draw(\i,9)node{$a_\i$};
\draw[draw=magenta,ultra thick] (1,4)to(1,5); \draw(1-0.3,5-0.3)node{$\magenta{x_1}$};
\draw[draw=magenta,ultra thick] (1,5)to(2,5);
\draw[draw=magenta,ultra thick] (2,5)to(2,6); \draw(2-0.3,6-0.3)node{$\magenta{\ov{x}_1}$};
\draw[draw=magenta,ultra thick] (2,6)to(3,6);
\draw[draw=magenta,ultra thick] (3,6)to(3,7); \draw(3-0.3,7-0.3)node{$\magenta{x_2}$};
\draw[draw=magenta,ultra thick] (3,7)to(8,7);
\draw[draw=magenta,ultra thick] (8,7)to(8,8); \draw(8-0.3,8-0.3)node{$\magenta{\ov{x}_4}$};
\draw[draw=blue,ultra thick] (3,3)to(6,3);
\draw[draw=blue,ultra thick] (6,3)to(6,4);  \draw(6-0.3,4-0.3)node{$\blue{\ov{x}_3}$};
\draw[draw=blue,ultra thick] (6,4)to(7,4);
\draw[draw=blue,ultra thick] (7,4)to(7,5);  \draw(7-0.3,5-0.3)node{$\blue{x_4}$};
\draw[draw=blue,ultra thick] (7,5)to(7,6);  \draw(7-0.3,6-0.3)node{$\blue{x_4}$};
\draw[draw=blue,ultra thick] (7,6)to(8,6);
\draw[draw=red,ultra thick] (5,2)to(7,2);
\draw[draw=red,ultra thick] (7,2)to(7,3);  \draw(7-0.3,3-0.3)node{$\red{x_4}$};
\draw[draw=red,ultra thick] (7,3)to(8,3);
\draw[draw=red,ultra thick] (8,3)to(8,4);  \draw(8-0.3,4-0.3)node{$\red{\ov{x}_4}$};
\draw[draw=red,ultra thick] (8,4)to(8,5);  \draw(8-0.3,5-0.3)node{$\red{\ov{x}_4}$};
\draw[draw=cyan,very thick] (7,1)to(8,1);
\foreach \j in {4,...,8} \draw(1,\j)node{$\bullet$};
\foreach \j in {4,...,8} \draw(2,\j)node{$\bullet$};
\foreach \j in {3,...,8} \draw(3,\j)node{$\bullet$};
\foreach \j in {3,...,8} \draw(4,\j)node{$\bullet$};
\foreach \j in {2,...,8} \draw(5,\j)node{$\bullet$};
\foreach \j in {2,...,8} \draw(6,\j)node{$\bullet$};
\foreach \j in {1,...,8} \draw(7,\j)node{$\bullet$};
\foreach \j in {1,...,8} \draw(8,\j)node{$\bullet$};
\end{tikzpicture}
}}\cr\cr
T\ = \
\YT{0.2in}{0.2in}{}{
 {\magenta{1},\magenta{\ov1},\magenta{2},\magenta{\ov4}},
 {\blue{\ov3},\blue{4},\blue{4}},
 {\red{4},\red{\ov4},\red{\ov4}},
}
\qquad\qquad
\wgt(T)\ = \ 
\YT{0.2in}{0.5in}{}{
 {\magenta{x_1},\magenta{\ov{x}_1},\magenta{{x_2}\!+\!{a_1}},\magenta{x_{\ov4}\!+\!{a_7}}},
 {\blue{x_{\ov3}\!+\!{a_1}},\blue{{x_4}\!+\!{a_3}},\blue{{x_4}\!+\!{a_4}}},
 {\red{{x_4}\!+\!{a_1}},\red{{x_{\ov4}}\!+\!{a_3}},\red{{x_{\ov4}}\!+\!{a_4}}},
} \nonumber
\end{array}
$$
\caption{Contribution to $sp_\lambda(\x,\ovx\,|\,\a)$ from $T$ and $LP(T)$, where $a_m=0$ for $m\leq0$.}
\label{fig-sp-TtoLP}
\end{figure*}

Finally, in the factorial odd orthogonal case the alphabet is extended to include not only both $x_k$ and $\ov{x}_k$ for $k=1,2,\ldots,n$, 
but also $1$, and $\a$ is replaced this time by $\tau^{1-n}\a$ as 
dictated by (\ref{eqn-so-hmtau}). The starting points are $P_i=(2i-1,n-i+1)$, 
ensuring as in the symplectic case that the condition {\bf(T4)} is satisfied, and the end points are $Q_j=(2n+1,n-j+1+\lambda_j)$ 
since the alphabet is now of length $2n+1$. To take into account the last factor $(1-a_{m+n})$
appearing in (\ref{eqn-so-hm-expansion}) the lattice paths may now include a final diagonal step.
The fact that it is diagonal ensures that there is at most one of these steps on each
lattice path. Once again it is only the $n$-tuples of non-intersecting lattice paths from $P_i$ to $Q_i$ that
contribute to $so_\lambda(\x,\ovx,1\,|\,\a)$ and these are in bijective correspondence with
the odd orthogonal tableaux of Definition~\ref{Def-soT} of shape $\lambda$ with entries from $\{1<\ov1<\cdots<n<\ov{n}<0\}$.
The fact that on each lattice path the final step is either vertical or diagonal, with the
latter to be associated with entries $0$ ensures that the condition {\bf(T5)} is automatically satisfied.
This is exemplified in Figure~\ref{fig-so-TtoLP} for $n=4$ and $\lambda=(4,3,3)$.
\qed

\begin{figure*}[h]  
$$
\begin{array}{c}
LP(T)\ =\ 
\vcenter{\hbox{
\begin{tikzpicture}[x={(0in,-0.3in)},y={(0.3in,0in)}] 
\foreach \i in {1,...,4} \draw[-](2*\i-1,5-\i)to(0,4+\i); \foreach \i in {1,...,3} \draw(0-0.2,8-\i+0.2)node{$a_{\ov{\i}}$}; 
\draw(0-0.2,8+0.2)node{$a_{0}$}; 
\foreach \i in {1,...,8} \draw[-](9-0.5,\i-0.5)to(\i-1,9); \foreach \i in {1,...,8} \draw(\i-1.2,9+0.2)node{$a_\i$};
\foreach \i in {1,...,8} \draw[-](9-0.5,\i-0.5)to(9+0.8,\i+0.8); \foreach \i in {1,...,8} \draw(10,\i+0.8)node{$-\!a_\i$};
\draw(1,4-0.4)node{$P_1$};
\draw(3,3-0.4)node{$P_2$};
\draw(5,2-0.4)node{$P_3$};
\draw(7,1-0.4)node{$P_4$};
\draw(9+0.4,1)node{$Q_4$};
\draw(9+0.4,5)node{$Q_3$};
\draw(9+0.4,6)node{$Q_2$};
\draw(9+0.4,8)node{$Q_1$};
\draw[draw=magenta,ultra thick] (1,4)to(1,5); \draw(1-0.3,5-0.3)node{$\magenta{x_1}$};
\draw[draw=magenta,ultra thick] (1,5)to(2,5);
\draw[draw=magenta,ultra thick] (2,5)to(2,6); \draw(2-0.3,6-0.3)node{$\magenta{\ov{x}_1}$};
\draw[draw=magenta,ultra thick] (2,6)to(3,6);
\draw[draw=magenta,ultra thick] (3,6)to(3,7); \draw(3-0.3,7-0.3)node{$\magenta{x_2}$};
\draw[draw=magenta,ultra thick] (3,7)to(8,7);
\draw[draw=magenta,ultra thick] (8,7)to(8,8); \draw(8-0.3,8-0.3)node{$\magenta{\ov{x}_4}$};
\draw[draw=magenta,ultra thick] (8,8)to(9,8);
\draw[draw=blue,ultra thick] (3,3)to(5,3);
\draw[draw=blue,ultra thick] (5,3)to(5,4);  \draw(5-0.3,4-0.3)node{$\blue{x_3}$};
\draw[draw=blue,ultra thick] (5,4)to(6,4);  
\draw[draw=blue,ultra thick] (6,4)to(7,4);
\draw[draw=blue,ultra thick] (7,4)to(7,5);  \draw(7-0.3,5-0.3)node{$\blue{x_4}$};
\draw[draw=blue,ultra thick] (7,5)to(8,5);  
\draw[draw=blue,ultra thick] (8,5)to(9,6);  \draw(9-0.4,6+0.0)node{$\blue{1}$};
\draw[draw=red,ultra thick] (5,2)to(7,2);
\draw[draw=red,ultra thick] (7,2)to(7,3);  \draw(7-0.3,3-0.3)node{$\red{x_4}$};
\draw[draw=red,ultra thick] (7,3)to(8,3);
\draw[draw=red,ultra thick] (8,3)to(8,4);  \draw(8-0.3,4-0.3)node{$\red{\ov{x}_4}$};
\draw[draw=red,ultra thick] (8,4)to(9,5);  \draw(9-0.4,5+0.0)node{$\red{1}$};
\draw[draw=cyan,ultra thick] (7,1)to(8,1);
\draw[draw=cyan,ultra thick] (8,1)to(9,1); 
\foreach \j in {4,...,8} \draw(1,\j)node{$\bullet$};
\foreach \j in {4,...,8} \draw(2,\j)node{$\bullet$};
\foreach \j in {3,...,8} \draw(3,\j)node{$\bullet$};
\foreach \j in {3,...,8} \draw(4,\j)node{$\bullet$};
\foreach \j in {2,...,8} \draw(5,\j)node{$\bullet$};
\foreach \j in {2,...,8} \draw(6,\j)node{$\bullet$};
\foreach \j in {1,...,8} \draw(7,\j)node{$\bullet$};
\foreach \j in {1,...,8} \draw(8,\j)node{$\bullet$};
\foreach \j in {1,...,8} \draw(9,\j)node{$\bullet$};
\end{tikzpicture}
}}\cr\cr
T\ =\ \YT{0.2in}{0.2in}{}{
 {\magenta{1},\magenta{\ov1},\magenta{2},\magenta{\ov4}},
 {\blue{3},\blue{4},\blue{0'}},
 {\red{4},\red{\ov4},\red{0'}},
}
\qquad\qquad
\wgt(T)\ =\ \YT{0.2in}{0.5in}{}{
 {\magenta{x_1},\magenta{\ov{x}_1},\magenta{{x_2}\!+\!{a_2}},\magenta{{x_{\ov4}}\!+\!{a_8}}},
 {\blue{{x_3}\!+\!{a_{1}}},\blue{{x_4}\!+\!{a_4}},\blue{{1}-{a_6}}},
 {\red{{x_4}\!+\!{a_2}},\red{{x_{\ov4}}\!+\!{a_4}},\red{{1}-{a_5}}},
}
$$ 
\end{array}
$$
\caption{Contribution to $so_\lambda(\x,\ovx\,|\,\a)$ from $T$ and $LP(T)$, where $a_m=0$ for $m\leq0$.}
\label{fig-so-TtoLP}
\end{figure*}

\section{Primed shifted tableaux and factorial $Q$-functions}\label{sec-tabP-factQ}

The passage from Schur functions to Schur $Q$-functions can be effected by replacing tableaux by primed shifted tableaux~\cite{Wor84,Sag87}.
We replicate this in the factorial setting by offering definitions of three types of factorial $Q$-functions expressed in terms of 
certain primed shifted tableaux. To this end we first define shifted Young diagrams.

A partition is said to be strict if its non-zero parts are distinct. 
Each such strict partition $\lambda$ of length $\ell(\lambda)\leq n$ specifies a shifted Young diagram $SF^\lambda$
consisting of rows of boxes of lengths $\lambda_i$ for $i=1,2,\ldots,\ell(\lambda)$ left adjusted to
a diagonal line. This is exemplified in the case $\lambda=(6,4,3)$ by
$$
SF^{6531}\ = \
\SYT{0.2in}{0.2in}{}{
 {{},{},{},{},{},{}},
 {{},{},{},{}},
 {{},{},{}},
}
$$
This allows us to define various primed shifted tableaux. 

\begin{Definition}\label{Def-gl-Qtab}~\cite{Wor84,Sag87} Let ${\cal P}^{\gl}_\lambda$ be the set of all primed shifted tableaux $P$
of shape $\lambda$ that are obtained by filling each box of $SF^\lambda$ with an entry $P_{ij}$
from the alphabet 
$$\{1'< 1< 2'<2<\cdots<n'<n\}$$
with one entry in each box, in such a way that:
{\bf (Q1)} entries weakly increase from left to right across rows;
{\bf (Q2)} entries weakly increase from top to bottom down columns;
{\bf (Q3)} no two identical unprimed entries appear in any column;
{\bf (Q4)} no two identical primed entries appear in any row;
\end{Definition}

\begin{Definition}\label{Def-sp-Qtab}~\cite{HK07} Let ${\cal P}^{\sp}_\lambda$ be the set of all primed shifted tableaux $P$
of shape $\lambda$ that are obtained by filling each box of $SF^\lambda$ with an entry $P_{ij}$
from the alphabet 
$$\{1'\!<1<\!\ov{1}'\!<\!\ov{1}\!<\!2'\!<\!2\!<\!\ov{2}'\!<\!\ov{2}<\cdots<n'<n<\ov{n}'<\ov{n}\}$$ 
with one entry in each box, 
in such a way that the conditions {\bf (Q1)-(Q4)} are satisfied together with:
{\bf (Q5)} at most one of $\{k',k,\ov{k}',\ov{k}\}$ appears on the main diagonal for each $k=1,2,\ldots,n$.
\end{Definition}

\begin{Definition}\label{Def-so-Qtab} Let ${\cal P}^{\so}_\lambda$ be the set of all primed shifted tableaux $P$
of shape $\lambda$ that are obtained by filling each box of $SF^\lambda$ with an entry $P_{ij}$
from the alphabet 
$$\{1'<1<\ov{1}'<\ov{1}<2'<2<\ov{2}'<\ov{2}<\cdots<n'<n<\ov{n}'<\ov{n}<0'\}$$
with one entry in each box, 
in such a way that the conditions {\bf (Q1)-(Q5)} are satisfied. together with:
{\bf (Q6)} the entry $0'$ does not appear on the main diagonal.
\end{Definition}

In the case $\lambda=(6,5,3)$ each of these types of shifted primed tableaux is illustrated by
\begin{equation}\label{eqn-Qtabx3}
\SYT{0.2in}{0.2in}{}{
 {{1'},{1},{2'},{2},{3'},{4}},
   {{2},{3'},{3},{3}},
     {{4'},{4},{4}},
}
\qquad
\SYT{0.2in}{0.2in}{}{
 {{1},{\ov1},{2'},{\ov2'},{3'},{3}},
   {{2'},{2},{3},{4'}},
     {{\ov4'},{4},{4}},
}
\qquad
\SYT{0.2in}{0.2in}{}{
 {{1},{\ov1},{2'},{\ov2'},{3},{0'}},
   {{\ov2'},{\ov2},{3},{4'}},
     {{4'},{4},{0'}},
}
\end{equation}

We then propose the following definitions of factorial $Q$-functions:
\begin{Definition}\label{Def-Q}
For $\a=(a_1,a_2,\ldots)$, $a_0=0$, and any strict partition $\lambda$ of length $\ell(\lambda)\leq n$, let 
\begin{equation}\label{eqn-def-Q}
   Q^\g_\lambda(\z;\w\,|\,\a) = \sum_{P\in {\cal P}^\g_\lambda}\ \prod_{(i,j)\in SF^\lambda}\!\! \wgt(P_{ij}) ~~~~~\mbox{where} 
\end{equation}
\begin{equation}\label{eqn-Pij-Q}
\begin{array}{|l|l|}
\hline
\g& Q^\g_\lambda(\z;\w\,|\,\a)\cr
\hline
gl&Q_\lambda(\x;\y\,|\,\a)\cr
\sp&Q^{\sp}_\lambda(\x,\ovx;\y,\ovy\,|\,\a)\cr
\so&Q^{\so}_\lambda(\x,\ovx;\y,\ovy,1\,|\,\a)\cr
\hline
\end{array}
~~~~~~~\mbox{and}~~~~~~~
\begin{array}{|l|l||l|l|}
\hline
P_{ij}&\wgt(P_{ij}) & P_{ij}&\wgt(P_{ij})\cr
\hline
k&x_k+a_{j-i} &  k'&y_k-a_{j-i}  \cr
 \ov{k}&\ov{x}_k+a_{j-i} &  \ov{k}'&\ov{y}_k-a_{j-i}\cr
& & 0'&1-a_{j-i} \cr
\hline
\end{array}
\end{equation}
\end{Definition}

It might be noted that for these factorial $Q$-functions the dependence on the factorial parameters $\a$ is simpler than it is for factorial characters
since the factors in (\ref{eqn-Pij-Q}) are all of the form $z_k\pm a_{j-i}$ with the subscript on $a$ 
completely independent of that on $z$.
  
The definition given here of $Q_\lambda(\x;\y\,|\,\a)$ has been introduced and studied elsewhere~\cite{HK15}. The special case 
$Q_\lambda(\x;\x\,|\,-\a)$ obtained by setting $\y=\x$ and $\a=-\a$ coincides with the generalized $Q$-function 
$Q_\lambda(\x\,|\,\a)$ introduced by Ivanov~\cite{Iva01,Iva05} and studied further by Ikeda, Milhalcea and Naruse~\cite{IMN11}.
If one further sets $\a=\0$ one recovers the combinatorial primed shifted tableaux formula~\cite{Wor84,Sag87,Ste89,Oka90}
for the original Schur $Q$-functions $Q_\lambda(\x)$.

\section{Determinantal expressions for factorial $Q$-functions}\label{sec-PtoLP-detQ}
In order to establish algebraic expressions for our factorial $Q$-functions we follow the method of Okada~\cite{Oka90} to 
construct lattice path models, in which each row of a primed shifted tableau, $P$, specifies a lattice path contributing to
an $\ell(\lambda)$-tuple, $LP(P)$, of non-intersecting lattice paths. 
In these models the $i$th paths extend from $P_i$ to $Q_i$ for $i=1,2,\ldots,\ell(\lambda)$ as specified by
\begin{equation}\label{eqn-PiQi}
\begin{array}{|c|c|c|c|}
\hline 
g&P_{ii}&P_i&Q_i\cr
\hline
gl&k,k'&(k,0)&(n,\lambda_i)\cr
\hline
\sp&k,k',\ov{k},\ov{k}'&(2k-\frac12,0)&(2n,\lambda_i)\cr
\hline
\so&k,k',\ov{k},\ov{k}'&(2k-\frac12,0)&(2n+1,\lambda_i)\cr
\hline
\end{array}
\end{equation}
It is convenient to set $d_i=k$ if $P_{ii}\in\{k,k',\ov{k},\ov{k}'\}$ and introduce ${\bf d}=(d_1,d_2,\ldots,d_{\ell(\lambda)})$. 

The entries $P_{ij}$ in the $i$th row of $P$ determine $\lambda_i$
edges of the corresponding lattice path extending from $P_i$ to $Q_i$.
The nature of these edges and their corresponding weights,
as determined by Definition~\ref{Def-Q}, are as prescribed below.

\begin{equation}\label{eqn-LPwgts}
\begin{array}{|c|c|c|c|c|c|c|}
\hline 
gl&P_{ii}=k&x_k&
\vcenter{\hbox{\begin{tikzpicture}[x={(0in,-0.25in)},y={(0.25in,0in)}] 
\draw[draw=black,ultra thick] (0,0)to[out=60,in=120](0,1); \draw(0,0)node{$\bullet$};\draw(0,1)node{$\bullet$};
\end{tikzpicture}}}
&P_{ii}=k'&y_k&
\vcenter{\hbox{\begin{tikzpicture}[x={(0in,-0.25in)},y={(0.25in,0in)}] 
\draw[draw=black,ultra thick] (0,0)to[out=-60,in=-120](0,1);  \draw(0,0)node{$\bullet$};\draw(0,1)node{$\bullet$};
\end{tikzpicture}}}
\cr
\hline\hline
\sp,\so&P_{ii}=k&x_k&
\vcenter{\hbox{\begin{tikzpicture}[x={(0in,-0.25in)},y={(0.25in,0in)}] 
\draw[draw=black,ultra thick] (0,0)to[out=45,in=180](0-0.5,1); \draw(0,0)node{$\bullet$};\draw(0-0.5,1)node{$\bullet$};
\end{tikzpicture}}}
&P_{ii}=k'&y_k&
\vcenter{\hbox{\begin{tikzpicture}[x={(0in,-0.25in)},y={(0.25in,0in)}] 
\draw[draw=black,ultra thick] (0,0)to[out=0,in=225](0-0.5,1);  \draw(0,0)node{$\bullet$};\draw(0-0.5,1)node{$\bullet$};
\end{tikzpicture}}}
\cr
\hline 
\sp,\so&P_{ii}=\ov{k}&\ov{x}_k&
\vcenter{\hbox{\begin{tikzpicture}[x={(0in,-0.25in)},y={(0.25in,0in)}] 
\draw[draw=black,ultra thick] (0,0)to[out=0,in=135](0+0.5,1); \draw(0,0)node{$\bullet$};\draw(0+0.5,1)node{$\bullet$};
\end{tikzpicture}}}
&P_{ii}=\ov{k}'&\ov{y}_k&
\vcenter{\hbox{\begin{tikzpicture}[x={(0in,-0.25in)},y={(0.25in,0in)}] 
\draw[draw=black,ultra thick] (0,0)to[out=315,in=180](0+0.5,1); \draw(0,0)node{$\bullet$};\draw(0+0.5,1)node{$\bullet$};
\end{tikzpicture}}}
\cr
\hline 
\sp,\so&P_{ij}=k, i<j& x_k+a_{j-i}&
\vcenter{\hbox{\begin{tikzpicture}[x={(0in,-0.25in)},y={(0.25in,0in)}] 
\draw[draw=black,ultra thick] (0,0)to(0,1); \draw(0,0)node{$\bullet$};\draw(0,1)node{$\bullet$};
\end{tikzpicture}}}
&P_{ij}=k', i<j&y_k-a_{j-i}&
\vcenter{\hbox{\begin{tikzpicture}[x={(0in,-0.25in)},y={(0.25in,0in)}] 
\draw[draw=black,ultra thick] (0,0)to(1,1); \draw(0,0)node{$\bullet$};\draw(1,1)node{$\bullet$};
\end{tikzpicture}}}
\cr
\hline
\sp,\so&P_{ij}=\ov{k}, i<j&\ov{x}_k+a_{j-i}&
\vcenter{\hbox{\begin{tikzpicture}[x={(0in,-0.25in)},y={(0.25in,0in)}] 
\draw[draw=black,ultra thick] (0,0)to(0,1); \draw(0,0)node{$\bullet$};\draw(0,1)node{$\bullet$};
\end{tikzpicture}}}
&P_{ij}=\ov{k}', i<j&\ov{y}_k-a_{j-i}&
\vcenter{\hbox{\begin{tikzpicture}[x={(0in,-0.25in)},y={(0.25in,0in)}] 
\draw[draw=black,ultra thick] (0,0)to(1,1); \draw(0,0)node{$\bullet$};\draw(1,1)node{$\bullet$};
\end{tikzpicture}}}
\cr
\hline\hline
\so&&&
&P_{ij}=0'&1-a_{j-i}&
\vcenter{\hbox{\begin{tikzpicture}[x={(0in,-0.25in)},y={(0.25in,0in)}] 
\draw[draw=black,ultra thick] (0,0)to(1,1); \draw(0,0)node{$\bullet$};\draw(1,1)node{$\bullet$};
\end{tikzpicture}}}
\cr
\hline 
\end{array}
\end{equation}
Each edge has its rightmost vertex at $(r,\ell)$ in the rectangular lattice,
with $\ell=0$ for the curved edges and $\ell=j-i$ for all the others,
while $r=k$ for entries $k$ and $k'$ in the case $\g=gl$, 
but $r=2k-1$ or $r=2k$ according as the entry $k$, $k'$, $\ov{k}$ or $\ov{k}'$ is 
unbarred or barred in the cases $\g=\sp$ and $\oo$, and $r=2n+1$
if the entry is $0'$ as may only occur in the case $\g=\oo$.
The path $P_iQ_i$ is completed by the insertion of vertical edges of weight $1$.
All this is illustrated in the case of our three running examples 
in Figures~\ref{fig-glQ-PtoLP}, \ref{fig-spQ-PtoLP}, and \ref{fig-soQ-PtoLP}, respectively, in the case $\lambda=(6,4,3)$,
for which $\ell(\lambda)=3$.

\begin{figure*}[htbp]
\begin{center}
$$
\begin{array}{c}
LP(P)\ =\ 
\vcenter{\hbox{
\begin{tikzpicture}[x={(0in,-0.3in)},y={(0.3in,0in)}] 
\draw(1,-0.4)node{$P_1$};
\draw(2,-0.4)node{$P_2$};
\draw(4,-0.4)node{$P_3$};
\draw(4.4,3)node{$Q_3$};
\draw(4.4,4)node{$Q_2$};
\draw(4.4,6)node{$Q_1$};
\foreach \j in {0,...,5} \draw[-](0-0.2,\j+1)to(4,\j+1);
\foreach \j in {1,...,5} \draw(-0.5,\j+1)node{$a_{\j}$};
\draw(-0.5,1)node{$a_{0}$};
\draw[draw=magenta,ultra thick] (1,0)to[out=-60,in=-120](1,1); \draw(1-0.0,1-0.5)node{${y_1}$};
\draw[draw=magenta,ultra thick] (1,1)to(1,2); \draw(1-0.3,2-0.5)node{$\magenta{x_1}$};
\draw[draw=magenta,ultra thick] (1,2)to(2,3); \draw(2-0.7,3-0.3)node{$\magenta{y_2}$};
\draw[draw=magenta,ultra thick] (2,3)to(2,4); \draw(2-0.3,4-0.5)node{$\magenta{x_2}$};
\draw[draw=magenta,ultra thick] (2,4)to(3,5); \draw(3-0.7,5-0.3)node{$\magenta{y_3}$};
\draw[draw=magenta,ultra thick] (3,5)to(4,5);
\draw[draw=magenta,ultra thick] (4,5)to(4,6); \draw(4-0.3,6-0.5)node{$\magenta{x_4}$};
\draw[draw=blue,ultra thick] (2,0)to[out=60,in=120](2,1); \draw(2+0.0,1-0.5)node{${x_2}$};
\draw[draw=blue,ultra thick] (2,1)to(3,2); \draw(3-0.7,2-0.3)node{$\blue{y_3}$};
\draw[draw=blue,ultra thick] (3,2)to(3,3); \draw(3-0.3,3-0.5)node{$\blue{x_3}$};
\draw[draw=blue,ultra thick] (3,3)to(3,4); \draw(3-0.3,4-0.5)node{$\blue{x_3}$};
\draw[draw=red,ultra thick] (4,0)to[out=-60,in=-120](4,1);  \draw(4-0.0,1-0.5)node{${y_4}$};
\draw[draw=red,ultra thick] (4,1)to(4,2); \draw(4-0.3,2-0.5)node{$\red{x_4}$};
\draw[draw=red,ultra thick] (4,2)to(4,3); \draw(4-0.3,3-0.5)node{$\red{x_4}$};
\foreach \i in {1,...,4} \foreach \j in {0,...,6} \draw(\i,\j)node{$\bullet$};
\end{tikzpicture}
}}\cr\cr
P\ =\  
\SYT{0.2in}{0.2in}{}{
 {\magenta{1'},\magenta{1},\magenta{2'},\magenta{2},\magenta{3'},\magenta{4}},
   {\blue{2},\blue{3'},\blue{3},\blue{3}},
     {\red{4'},\red{4},\red{4}},
}
\quad
\wgt(P)\ = \ 
\SYT{0.2in}{0.5in}{}{
 {\magenta{y_1},\magenta{x_1\!+\!a_1},\magenta{y_2\!-\!a_2},\magenta{x_2\!+\!a_3},\magenta{y_3\!-\!a_4},\magenta{x_4\!+\!a_5}},
   {\blue{x_2},\blue{y_3\!-\!a_1},\blue{x_3\!+\!a_2},\blue{x_3\!+\!a_3}},
     {\red{y_4},\red{x_4\!+\!a_1},\red{x_4\!+\!a_2}},
}
\end{array}
$$
\end{center}
\caption{Contribution to $Q_\lambda(\x;\y\,|\,\a)$ from $P$ and $LP(P)$, where $a_0=0$.}
\label{fig-glQ-PtoLP}
\end{figure*}
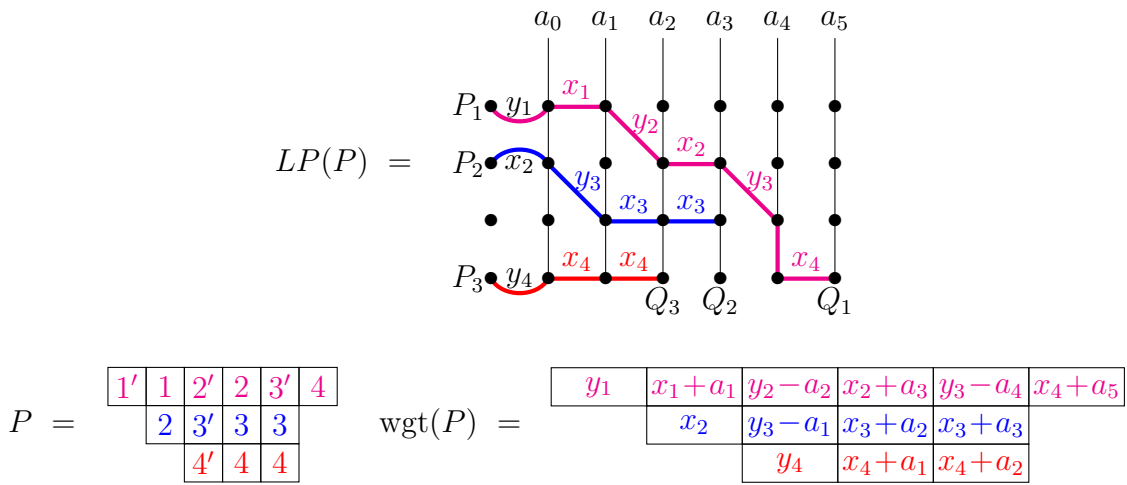

\begin{figure*}[htbp]
\begin{center}
$$
\begin{array}{c}
LP(P)\ =\ 
\vcenter{\hbox{
\begin{tikzpicture}[x={(0in,-0.3in)},y={(0.3in,0in)}] 
\draw(1.5,-0.4)node{$P_1$};
\draw(3.5,-0.4)node{$P_2$};
\draw(7.5,-0.4)node{$P_3$};
\draw(8.4,3)node{$Q_3$};
\draw(8.4,4)node{$Q_2$};
\draw(8.4,6)node{$Q_1$};
\foreach \j in {0,...,5} \draw[-](1-0.5,\j+1)to(8,\j+1);
\foreach \j in {1,...,5} \draw(1-0.7,\j+1)node{$a_{\j}$};
\draw(1-0.7,1)node{$a_0$};
\draw[draw=magenta,ultra thick] (1.5,0)to[out=45,in=180](1,1); \draw(1-0.1,1-0.6)node{$\magenta{x_1}$};
\draw[draw=magenta,ultra thick] (1,1)to(2,1);
\draw[draw=magenta,ultra thick] (2,1)to(2,2); \draw(2-0.3,2-0.5)node{$\magenta{\ov{x}_1}$};
\draw[draw=magenta,ultra thick] (2,2)to(3,3); \draw(3-0.7,3-0.3)node{$\magenta{y_2}$};
\draw[draw=magenta,ultra thick] (3,3)to(4,4); \draw(4-0.7,4-0.3)node{$\magenta{\ov{y}_2}$};
\draw[draw=magenta,ultra thick] (4,4)to(5,4);
\draw[draw=magenta,ultra thick] (5,4)to(5,5); \draw(5-0.3,5-0.5)node{$\magenta{x_3}$};
\draw[draw=magenta,ultra thick] (5,5)to(5,6); \draw(5-0.3,6-0.5)node{$\magenta{x_3}$};
\draw[draw=magenta,ultra thick] (5,6)to(8,6);
\draw[draw=blue,ultra thick] (3.5,0)to[out=0,in=135](4,1); \draw(4+0.0,1-0.5)node{$\blue{\ov{x}_2}$};
\draw[draw=blue,ultra thick] (4,1)to(4,2); \draw(4-0.3,2-0.5)node{$\blue{\ov{x}_2}$};
\draw[draw=blue,ultra thick] (4,2)to(5,2); 
\draw[draw=blue,ultra thick] (5,2)to(5,3); \draw(5-0.3,3-0.5)node{$\blue{x_3}$};
\draw[draw=blue,ultra thick] (5,3)to(6,3);
\draw[draw=blue,ultra thick] (6,3)to(7,4); \draw(7-0.7,4-0.3)node{$\blue{y_4}$};
\draw[draw=blue,ultra thick] (7,4)to(8,4);
\draw[draw=red,ultra thick] (7.5,0)to[out=0,in=225](7,1);  \draw(7-0.0,1-0.5)node{$\red{y_4}$};
\draw[draw=red,ultra thick] (7,1)to(7,2); \draw(7-0.3,2-0.5)node{$\red{x_4}$};
\draw[draw=red,ultra thick] (7,2)to(8,2); 
\draw[draw=red,ultra thick] (8,2)to(8,3); \draw(8-0.3,3-0.5)node{$\red{\ov{x}_4}$};
\foreach \i in {1,...,8} \foreach \j in {1,...,6} \draw(\i,\j)node{$\bullet$};
\foreach \i in {1,3,5,7} \draw(\i+0.5,0)node{$\bullet$};
\end{tikzpicture}
}}\cr\cr
P\ =\  
\SYT{0.2in}{0.2in}{}{
 {\magenta{1},\magenta{\ov1},\magenta{2'},\magenta{\ov2'},\magenta{3},\magenta{3}},
   {\blue{\ov2},\blue{\ov2},\blue{3},\blue{4'}},
     {\red{4'},\red{4},\red{\ov4}},
}
\quad
\wgt(P)\ = \ 
\SYT{0.2in}{0.5in}{}{
 {\magenta{x_1},\magenta{\ov{x}_1\!+\!a_1},\magenta{y_2\!-\!a_2},\magenta{\ov{y}_2\!-\!a_3},\magenta{x_3\!+\!a_4},\magenta{x_3\!+\!a_5}},
   {\blue{\ov{x}_2},\blue{\ov{x}_2\!+\!a_1},\blue{x_3\!+\!a_2},\blue{y_4\!-\!a_3}},
     {\red{y_4},\red{x_4\!+\!a_1},\red{\ov{x}_4\!+\!a_2}},
}
\end{array}
$$
\end{center}
\caption{Contribution to $Q^{\sp}_\lambda(\x,\ovx;\y,\ovy\,|\,\a)$ from $P$ and $LP(P)$, where $a_0=0$.}
\label{fig-spQ-PtoLP}
\end{figure*}

\begin{figure*}[htbp]
\begin{center}
$$
\begin{array}{c}
LP(P)\ =\ 
\vcenter{\hbox{
\begin{tikzpicture}[x={(0in,-0.3in)},y={(0.3in,0in)}] 
\draw(1.5,-0.4)node{$P_1$};
\draw(3.5,-0.4)node{$P_2$};
\draw(7.5,-0.4)node{$P_3$};
\draw(9.4,3)node{$Q_3$};
\draw(9.4,4)node{$Q_2$};
\draw(9.4,6)node{$Q_1$};
\foreach \j in {0,...,5} \draw[-](1-0.5,\j+1)to(9,\j+1);
\foreach \j in {1,...,5} \draw(1-0.7,\j+1)node{$a_{\j}$};
\draw(1-0.7,1)node{$a_0$};
\draw[draw=magenta,ultra thick] (1.5,0)to[out=45,in=180](1,1); \draw(1-0.1,1-0.6)node{$\magenta{x_1}$};
\draw[draw=magenta,ultra thick] (1,1)to(2,1);
\draw[draw=magenta,ultra thick] (2,1)to(2,2); \draw(2-0.3,2-0.5)node{$\magenta{\ov{x}_1}$};
\draw[draw=magenta,ultra thick] (2,2)to(3,3); \draw(3-0.7,3-0.3)node{$\magenta{y_2}$};
\draw[draw=magenta,ultra thick] (3,3)to(4,4); \draw(4-0.7,4-0.3)node{$\magenta{\ov{y}_2}$};
\draw[draw=magenta,ultra thick] (4,4)to(5,4);
\draw[draw=magenta,ultra thick] (5,4)to(5,5); \draw(5-0.3,5-0.5)node{$\magenta{x_3}$};
\draw[draw=magenta,ultra thick] (5,5)to(8,5); 
\draw[draw=magenta,ultra thick] (8,5)to(9,6); \draw(9-0.7,6-0.3)node{$\magenta{1}$};
\draw[draw=blue,ultra thick] (3.5,0)to[out=315,in=180](4,1); \draw(4+0.1,1-0.5)node{$\blue{\ov{y}_2}$};
\draw[draw=blue,ultra thick] (4,1)to(4,2); \draw(4-0.3,2-0.5)node{$\blue{\ov{x}_2}$};
\draw[draw=blue,ultra thick] (4,2)to(5,2); 
\draw[draw=blue,ultra thick] (5,2)to(5,3); \draw(5-0.3,3-0.5)node{$\blue{x_3}$};
\draw[draw=blue,ultra thick] (5,3)to(6,3);
\draw[draw=blue,ultra thick] (6,3)to(7,4); \draw(7-0.7,4-0.3)node{$\blue{y_4}$};
\draw[draw=blue,ultra thick] (7,4)to(9,4);
\draw[draw=red,ultra thick] (7.5,0)to[out=0,in=225](7,1);  \draw(7-0.0,1-0.5)node{$\red{y_4}$};
\draw[draw=red,ultra thick] (7,1)to(7,2); \draw(7-0.3,2-0.5)node{$\red{x_4}$};
\draw[draw=red,ultra thick] (7,2)to(8,2); 
\draw[draw=red,ultra thick] (8,2)to(9,3); \draw(9-0.7,3-0.3)node{$\red{1}$};
\foreach \i in {1,...,9} \foreach \j in {1,...,6} \draw(\i,\j)node{$\bullet$};
\foreach \i in {1,3,5,7} \draw(\i+0.5,0)node{$\bullet$};
\end{tikzpicture}
}}\cr\cr
P\ =\  
\SYT{0.2in}{0.2in}{}{
 {\magenta{1},\magenta{\ov1},\magenta{2'},\magenta{\ov2'},\magenta{3},\magenta{0'}},
   {\blue{\ov2'},\blue{\ov2},\blue{3},\blue{4'}},
     {\red{4'},\red{4},\red{0'}},
}
\quad
\wgt(P)\ = \ 
\SYT{0.2in}{0.5in}{}{
 {\magenta{x_1},\magenta{\ov{x}_1\!+\!a_1},\magenta{y_2\!-\!a_2},\magenta{\ov{y}_2\!-\!a_3},\magenta{x_3\!+\!a_4},\magenta{1\!-\!a_5}},
   {\blue{\ov{y}_2},\blue{\ov{x}_2\!+\!a_1},\blue{x_3\!+\!a_2},\blue{y_4\!-\!a_3}},
     {\red{y_4},\red{x_4\!+\!a_1},\red{1\!-\!a_2}},
}
\end{array}
$$
\end{center}
\caption{Contribution to $Q^{\oo}_\lambda(\x,\ovx;\y,\ovy,1\,|\,\a)$ from $P$ and $LP(P)$, where $a_0=0$.}
\label{fig-soQ-PtoLP}
\end{figure*}

The sets of all non-intersecting $\ell(\lambda)$-tuples of lattice paths from fixed starting points 
$P_i$ to fixed end points $Q_i$ determined by ${\bf d}$ and $\lambda$ as tabulated in (\ref{eqn-PiQi}) 
are in bijective correspondence with sets of all the primed shifted tableaux of Definitions~\ref{Def-gl-Qtab}-\ref{Def-so-Qtab}
with diagonal entries $P_{ii}$ consistent with $d_i$ for all $i=1,2,\ldots,d_{\ell(\lambda)}$.
The sets of $\ell(\lambda)$-tuples of lattice paths may be extended to include those from $P_i$ to $Q_{\pi(i)}$ for
$i=1,2,\ldots,\ell(\lambda)$ with $\pi$ any permutation in $S_{\ell(\lambda)}$. Assigning weights in
accordance with (\ref{eqn-LPwgts}), the sum of the contributions from all these $\ell(\lambda)$-tuples 
multiplied by $\sgn(\pi)$ then constitutes a determinant. As established by Okada~\cite{Oka90}
all contributions mutually cancel, other than those from the non-intersecting $\ell(\lambda)$-tuples
for which $\pi$ is the identity element. It then follows that our $Q$-functions may be expressed in the form
\begin{align}
   Q_\lambda(\x;\y\,|\,\a)&= \sum_{\bf d}\left|\,(x_{d_i}\!+\!y_{d_i})\
	       \tilde{q}_{\lambda_j-1}(\x^{(d_i)};\y^{(d_i+1)}\,|\,\a)\,\right|; \label{eqn-glQdet} \\
	 Q^{\sp}_\lambda(\x,\!\ovx;\y,\!\ovy\,|\,\a)&=\sum_{\bf d}  
	       \big|\,(x_{d_i}\!+\!y_{d_i})\ \tilde{q}_{\lambda_j-1}(\x^{(d_i)},\,\ovx^{(d_i)};\y^{(d_i+1)},\,\ovy^{(d_i)},|\,\a) \cr
				&{\phantom{\sum_{\bf d}}}+ (\ov{x}_{d_i}\!+\!\ov{y}_{d_i})\ \tilde{q}_{\lambda_j-1}(\x^{(d_i+1)},\,\ovx^{(d_i)};\y^{(d_i+1)},\,\ovy^{(d_i+1)}\,|\,\a)\,\big|; \label{eqn-spQdet} \\
	 Q^{\so}_\lambda(\x,\!\ovx;\y,\!\ovy\,1\,|\,\a)&=\sum_{\bf d}  
	      \big|\,(x_{d_i}\!+\!y_{d_i})\ \tilde{q}_{\lambda_j-1}(\x^{(d_i)},\,\ovx^{(d_i)};\y^{(d_i+1)},\,\ovy^{(d_i)},\,1\,|\,\a) \cr
				&{\phantom{\sum_{\bf d}}}+ (\ov{x}_{d_i}\!+\!\ov{y}_{d_i})\ \tilde{q}_{\lambda_j-1}(\x^{(d_i+1)},\,\ovx^{(d_i)};\y^{(d_i+1)},\,\ovy^{(d_i+1)},\,1\,|\,\a)\,\big|, \label{eqn-soQdet}
\end{align}
where for any $m\geq0$, and all relevant $\u=(u_1,u_2,\ldots,u_r)$ and $\v=(v_1,v_2,\ldots,v_s)$
\begin{equation}\label{eqn-qmuv}
     \tilde{q}_{m}(\u;\v\,|\,\a)=\sum_{1\dot\leq i_1\dot\leq i_2\dot\leq \cdots\dot\leq i_m\dot\leq r+s}\ 
		             (w_{i_1}\!\pm\!a_1)(w_{i_2}\!\pm\!a_2)\cdots(w_{i_m}\!\pm\!a_m)\,,
\end{equation}
with $\w=(w_1,w_2,\ldots,w_{r+s})$ to be identified with an alternating ordered sequence of all the elements of $\u$ and $\v$.
The notation $\dot\leq$ is intended to indicate that the summation allows factors $(w_i\pm a_\ell) = (u_k + a_\ell)$ or $(v_k - a_\ell)$ to appear
according as $w_i = u_k$ or $v_k$, with several factors of the form $(u_k + a_\ell)(u_k + a_{\ell+1})\cdots$
allowed, but at most one factor $(v_k-a_{\ell})$. 

To present this result in a neater form it is convenient to make use of the following
\begin{Definition} 
\begin{align}
   q^{\gl}_m(\x^{(d)};\y^{(d+1)}\,|\,\a) 
	&= [t^m]\ \frac{\prod_{j=d+1}^n (1\!+\!ty_j)\prod_{k=1}^{m}(1\!+\!ta_k)}{\prod_{i=d}^n (1\!-\!tx_i)};     \label{eqn-gl-qmd} \\
	 q^{\sp}_{m}(\x^{(d)},\!\ovx^{(d)};\y^{(d+1)},\!\ovy^{(d+1)}\,|\,\a)  
	&=  [t^m]\ \frac{\prod_{j=d+1}^n ((1\!+\!ty_j)(1\!+\!t\ov{y}_j))\prod_{k=1}^{m}(1\!+\!ta_k)}{\prod_{i=d}^n ((1\!-\!tx_i)(1\!-\!t\ov{x}_i))};  \label{eqn-sp-qmd}  \\
	 q^{\so}_{m}(\x^{(d)},\!\ovx^{(d)};\y^{(d+1)},\!\ovy^{(d+1)},\!1\,|\,\a) 
	&=  [t^m]\ \frac{(1\!+\!t)\prod_{j=d+1}^n ((1\!+\!ty_j)(1\!+\!t\ov{y}_j))\prod_{k=1}^{m}(1\!+\!ta_k)}{\prod_{i=d}^n ((1\!-\!tx_i)(1\!-\!t\ov{x}_i))}.\label{eqn-so-qmd} 
\end{align}
\end{Definition}
In terms of these we have
\begin{Theorem}\label{The-Q}
Let $\x=(x_1,x_2,\ldots,x_n)$ and  $\y=(y_1,y_2,\ldots,y_n)$. Then for
any strict partition $\lambda$ of length $\ell(\lambda)\leq n$ and any $\a=(a_1,a_2,\ldots)$
\begin{align}
   Q_\lambda(\x;\y\,|\,\a)&= \sum_{\bf d}\left|\,(x_{d_{i}}\!+\!y_{d_{i}})
	q^{\gl}_{\lambda_j-1}(\x^{(d_{i})};\y^{(d_{i}+1)}\,|\,\a)\,\right|; \label{eqn-glQ} \\
	 Q^{\sp}_\lambda(\x,\!\ovx;\y,\!\ovy\,|\,\a)&=\sum_{\bf d}  \left|\,(x_{d_{i}}\!+\!y_{d_{i}}\!+\!\ov{x}_{d_{i}}\!+\!\ov{y}_{d_{i}})
	q^{\sp}_{\lambda_j-1}(\x^{(d_{i})},\!\ovx^{(d_{i})};\y^{(d_{i}+1)},\!\ovy^{(d_{i}+1)}\,|\,\a)\,\right|;\label{eqn-spQ}\\
   Q^{\oo}_\lambda(\x,\!\ovx;\y,\!\ovy,\!1\,|\,\a)&= \sum_{\bf d} \left|\,(x_{d_{i}}\!+\!y_{d_{i}}\!+\!\ov{x}_{d_{i}}\!+\!\ov{y}_{d_{i}})
  q^{\oo}_{\lambda_j-1}(\x^{(d_{i})},\!\ovx^{(d_{i})};\y^{(d_{i}+1)},\!\ovy^{(d_{i}+1)},\!1\,|\,\a)\,\right|, \label{eqn-soQ}
\end{align}
where each determinant is $\ell(\lambda)\times\ell(\lambda)$ and each sum is over all ${\bf d}=(d_1,d_2,\ldots,d_{\ell(\lambda)})$ 
such that $1\leq d_1<d_2<\cdots<d_{\ell(\lambda)}\leq n$.
\end{Theorem}

\noindent{\bf Proof}: It should first be noted that in the $\a=\0$ case it is known~\cite{Oka90} that for all $\u=(u_1,u_2,\ldots,u_r)$
and $\v=(v_1,v_2,\ldots,v_s)$ that
\begin{equation}
\tilde{q}_{m}(\u;\v\,|\,\0) = [t^m]\ \prod_{i=1}^r \frac{1}{1-tu_i}\  \prod_{j=1}^s\,(1+tv_j)\,.
\end{equation} 
More generally, consider 
\begin{equation}\label{eqn-qmuv-gf}
\tilde{q}_{m}(\u;\v\,|\,\a) = [t^m]\ \prod_{i=1}^r \frac{1}{1-tu_i}\ \prod_{j=1}^s\,(1+tv_j)\  \prod_{k=1}^{m+r-s-1} (1+ta_k)\,.
\end{equation}
By writing $(1+ta_{m+r-s-1})/(1-tu_r)=1+t(u_r+a_{m+r-s-1})/(1-tu_r)$ and $(1+tv_s)=(1+ta_{m+r-s})+t(v_s-a_{m+r-s})$ 
it may be verified, very much as in (\ref{eqn-hm-recur}) and (\ref{eqn-so-hmtau}), that 
\begin{equation}\label{eqn-qm-ur}
\tilde{q}_{m}(\u;\v\,|\,\a) = \tilde{q}_{m}(\u';\v\,|\,\a)+(u_r+a_{m+r-s-1})\tilde{q}_{m-1}(\u;\v\,|\,\a)
\end{equation}
and
\begin{equation}\label{eqn-qm-vs}
\tilde{q}_{m}(\u;\v\,|\,\a) = \tilde{q}_{m}(\u;\v'\,|\,\a)+(v_s-a_{m+r-s})\tilde{q}_{m-1}(\u;\v'\,|\,\a)
\end{equation}
where $\u'=(u_1,u_2,\ldots,u_{r-1})$ and $\v'=(v_1,v_2,\ldots,v_{s-1})$.
Applying the first of these in the case $\u=\x^{(d)}$ and $\v=\y^{(d+1)}$, for which $r=s+1$, leads to a factor $(x_n+a_m)$
in the right hand term.
Repeating this process for $\tilde{q}_{m-1}(\u;\v\,|\,\a)$ but still with $\u=\x^{(d)}$ and $\v=\y^{(d+1)}$, 
leads to a further factor $(x_n+a_{m-1})$ and so. In this way any dependence on $x_n$ takes the
form $(x_n+a_{\ell+1})\cdots(x_n+a_{m-1})(x_n+a_m)$. 
Then applying (\ref{eqn-qm-vs}) in the case $m=\ell$, $\u=\x'^{(d)}$ and $\v=\y^{(d+1)}$, for which $r=s$, leads to a factor $(y_n+a_\ell)$,
with no further factors involving $y_n$ allowed. 
Iterating these results it can be seen that $q^{\gl}_m(\x^{(d)};\y^{(d+1)}\,|\,a)$,
as defined in (\ref{eqn-gl-qmd}), generates the expansion (\ref{eqn-qmuv}) of $\tilde{q}_{m}(\u;\v\,|\,\a)$ in the case $\u=\x^{(d)}$
and $\v=\y^{(d+1)}$, thereby proving (\ref{eqn-glQ}). 
Moreover, in the symplectic case it can be shown directly from the generating functions appearing in 
(\ref{eqn-qmuv-gf}) and (\ref{eqn-sp-qmd}) that 
\begin{align}\label{eqn-qq-qsp}
  &(x_i+y_i)\ \tilde{q}_{\lambda_j-1}(x^{(i)},\ov{x}^{(i)};y^{(i+1)},\ov{y}^{(i)} \,|\, \a)
  + (\ov{x}_i +\ov{y}_i)\ \tilde{q}_{\lambda_j-1}(x^{(i+1)},\ov{x}^{(i)};y^{(i+1)},\ov{y}^{(i+1)} \,|\, \a)\cr
	&~~~~= (x_i+y_i+\ov{x}_i +\ov{y}_i)\  q^{\sp}_{\lambda_j-1}(x^{(i)},\ov{x}^{(i)};y^{(i+1)},\ov{y}^{(i+1)} \,|\, \a)\,,
\end{align}
thereby proving (\ref{eqn-spQ}), with a similar result applying to the odd orthogonal case.
\qed

It might be noted that each of the expressions in Theorem~\ref{The-Q} in the form of a sum over determinants
may be expressed directly as Pfaffian following, for example, the prescription for dealing with non-intersecting
lattice paths from a selection of fixed starting points to fixed set of end points~\cite{Ste90}. This has been done 
already in the case of the factorial $Q$-functions~\cite{Iva05,IMN11}, and will not be pursued here, but will be the subject of future work.
 
\section{Tokuyama identities}\label{sec-Tok}
Here we restrict ourselves to the case for which $\lambda=\mu+\delta$ where $\mu$ is a partition
of length $\ell(\mu)\leq n$ and $\delta=(n,n-1,\ldots,1)$ so that $\lambda$ is a strict partition of length
$\ell(\lambda)=n$. In such case the sums over ${\bf d}$ appearing in Theorem~\ref{The-Q}
reduce to a single term corresponding to the only possible case ${\bf d}=(1,2,\ldots,n)$. Moreover,
each of the surviving determinants factorises, to yield the following factorial Tokuyama type identities.


\begin{Theorem}\label{The-fTok}
Let $\lambda=\mu+\delta$ with $\delta=(n,n-1,\ldots,1)$ and $\mu$ a partition of length~~$\ell(\mu)\leq n$. 
Then for any $\x=(x_1,x_2,\ldots,x_n)$, $\y=(y_1,y_2,\ldots,y_n)$ and $\a=(a_1,a_2,\ldots)$
\begin{align}
    Q^{\gl}_\lambda(\x;\y\,|\,\a) &= \prod_{1\leq i\leq j\leq n} (x_i+y_j)\ s_\mu(\x\,|\,\a)\,; \label{eqn-gl-fTok} \\
		Q^{\sp}_\lambda(\x,\ov{\x};\y,\ov{\y}\,|\,\a) &= \prod_{1\leq i\leq j\leq n} (x_i+y_j+\ov{x}_i+\ov{y}_j)\ sp_\mu(\x,\ov{\x}\,|\,\a)\,;  \label{eqn-sp-fTok} \\
    Q^{\oo}_\lambda(\x,\ov{\x};\y,\ov{\y},1\,|\,\a) &= \prod_{1\leq i\leq j\leq n} (x_i+y_j+\ov{x}_i+\ov{y}_j)\ so_\mu(\x,\ov{\x},1\,|\,\a)\,. \label{eqn-so-fTok} 
\end{align}
\end{Theorem}

Before embarking on the proof it is helpful to make the following definition:
\begin{Definition}\label{Def-fmpqn} For all $1\leq p \leq q \leq n$ let
\begin{align}
   f^{\gl}_{m,p,q,n}(\x;\y\,|\,\a) 
	&= [t^m]\ \frac{\prod_{j=q+1}^n (1+ty_j)\prod_{k=1}^{m+q-p} (1+ta_k)}{\prod_{i=p}^n (1-tx_i)} \ ; \label{eqn-gl-fmpqn} \\
	 f^{\sp}_{m,p,q,n}(\x,\ovx;\y,\ovy\,|\,\a) 
	&= [t^m]\ \frac{\prod_{j=q+1}^n ((1+ty_j)(1-t\ov{y}_j))\prod_{k=1}^{m+q-p} (1+ta_k)}{\prod_{i=p}^n ((1-tx_i)(1-t\ov{x}_i))} \ \,; \label{eqn-sp-fmpqn} \\
		 f^{\so}_{m,p,q,n}(\x,\ovx;\y,\ovy,1\,|\,\a) 
	&= [t^m]\ \frac{(1+t)\prod_{j=q+1}^n ((1+ty_j)(1-t\ov{y}_j))\prod_{k=1}^{m+q-p} (1+ta_k)}{\prod_{i=p}^n ((1-tx_i)(1-t\ov{x}_i))} \ \,. \label{eqn-so-fmpqn}
\end{align}
\end{Definition}

In the special case $p=q=d$ these definitions are such that
\begin{align}
   f^{\gl}_{m,d,d,n}(\x;\y\,|\,\a)&=q^{\gl}_m(\x^{(d)};\y^{(d+1)}\,|\,\a) \,; \label{eqn-gl-fddn} \\
	 f^{\sp}_{m,d,d,n}(\x,\ovx;\y,\ovy\,|\,\a)&= q^{\sp}_{m}(\x^{(d)},\!\ovx^{(d)};\y^{(d+1)},\!\ovy^{(d+1)}\,|\,\a) \,;\label{eqn-sp-fddn} \\
	 f^{\so}_{m,d,d,n}(\x,\ovx;\y,\ovy,1\,|\,\a) &= q^{\so}_{m}(\x^{(d)},\!\ovx^{(d)};\y^{(d+1)},\!\ovy^{(d+1)},\!1\,|\,\a) ;\label{eqn-so-fddn}
\end{align}
and in the case $p=d,q=n$ they reduce to
\begin{align}
   f^{\gl}_{m,d,n,n}(\x;\y\,|\,\a)&=h_m^{\gl}(\x^{(d)}\,|\,\a)\,; \label{eqn-gl-fdnn} \\
	 f^{\sp}_{m,d,n,n}(\x,\ovx;\y,\ovy\,|\,\a)&=h_m^{\sp}(\x^{(d)},\ovx^{(d)}\,|\,\a)\,;\label{eqn-sp-fdnn} \\
	 f^{\so}_{m,d,n,n}(\x,\ovx;\y,\ovy,1\,|\,\a) &=	h_m^{\so}(\x^{(d)},\ovx^{(d)},1\,|\,\a)\,;\label{eqn-so-fdnn}
\end{align}

Finally, for $1\leq p<q\leq n$
\begin{align}
    &f^{\gl}_{m,p,q-1,n}(\x;\y\,|\,\a)- f^{\gl}_{m,p+1,q,n}(\x;\y\,|\,\a)=(x_p+y_q) f^{\gl}_{m-1,p,q,n}(\x;\y\,|\,\a)\,; \label{eqn-gl-fdiff} \\
		&f^{\sp}_{m,p,q-1,n}(\x,\ovx;\y,\ovy\,|\,\a) - f^{\sp}_{m,p+1,q,n}(\x,\ovx;\y,\ovy\,|\,\a) \cr
		&~~~~~~~~~~~~~~= (x_p+y_q+\ov{x}_p + \ov{y}_q) f^{\sp}_{m-1,p,q,n}(\x,\ovx;\y,\ovy\,|\,\a)\,; \label{eqn-sp-fdiff} \\
		&f^{\so}_{m,p,q-1,n}(\x,\ovx;\y,\ovy,1\,|\,\a) -  f^{\so}_{m,p+1,q,n}(\x,\ovx;\y,\ovy,1\,|\,\a) \cr
		&~~~~~~~~~~~~~~= (x_p+y_q+\ov{x}_p + \ov{y}_q) f^{\so}_{m-1,p,q,n}(\x,\ovx;\y,\ovy,1\,|\,\a) \,. \label{eqn-so-fdiff}
\end{align}

\noindent{\bf Proof of Theorem~\ref{The-fTok}}: 
If we now focus, for example, on the symplectic case (\ref{eqn-sp-fTok}) 
and start by using (\ref{eqn-spQ}) in the case $\ell(\lambda)=n$ then, as we have said, 
the sum over ${\bf d}$ is restricted to a single term with $d_i=i$ for $i=1,2,\ldots,n$. It follows that
\begin{equation}
  Q^{\sp}_\lambda(\x,\ov{\x};\y,\ov{\y}\,|\,\a) =  \prod_{i=1}^n\, {(x_i\!+\!y_i\!+\!\ov{x}_i\!+\!\ov{y}_i)}
	                 \left|\, f^{\sp}_{\lambda_j-1;i,i,n}(\x,\ov{\x};\y,\ov{\y}\,|\a) \,\right|  \,,
\end{equation}
where we have extracted a common factor $(x_i+y_i+\ov{x}_i + \ov{y}_i)$ from the $i$th row for $i=1,2,\ldots,n$, and used (\ref{eqn-sp-fddn}).
Then, by the repeated subtraction of successive rows from one another and using (\ref{eqn-sp-fdiff}), precisely as in (\ref{eqn-gl-hlambda}), we have 
\begin{equation}
  \left| f^{\sp}_{\lambda_j-1;i,i,n}(\x,\ov{\x};\y,\ov{\y}\,|\a) \,\right| 
   = \prod_{1\leq i<j\leq n}\, {(x_i\!+\!y_j\!+\!\ov{x}_i\!+\!\ov{y}_j)}\ \left|\ f^{\sp}_{\lambda_j-1-n+i;i,n,n}(\x,\ov{\x};\y,\ov{\y}\,|\a) \,\right| \,.
\end{equation}
We are now in a position to use (\ref{eqn-sp-fdnn}) which leads directly to 
\begin{align}	
 &f^{\sp}_{\lambda_j-1-n+i;i,n,n}(\x,\ov{\x};\y,\ov{\y}\,|\a) 
 = \left|\, h^{\sp}_{\lambda_j-(n-i+1)}(\x^{(i)},\ovx^{(i)}\,|\,\a) \,\right|  \cr
 &=\left|\, h^{\sp}_{\mu_j-j+i}(\x^{(i)},\ovx^{(i)}\,|\,\a) \,\right| 
 = sp_{\mu}(\x,\ov{\x}\,|\,\a) \,,																			
\end{align}
as required to complete the proof of (\ref{eqn-sp-fTok}). 
The final steps exploit the fact that $\lambda_j=\mu_j+n-j+1$ for $j=1,2,\ldots,n$,
as well as the symplectic factorial flagged Jacobi-Trudi identity of Theorem~\ref{The-fJT}. 
The other results (\ref{eqn-gl-fTok}) and (\ref{eqn-so-fTok}) can be established in exactly the same way.
\qed

\noindent{\bf Acknowledgements}
\label{sec:ack}
This work was supported by the Canadian Tri-Council Research
Support Fund. The first author (AMH) acknowledges the support of a
Discovery Grant from the Natural Sciences and Engineering Research Council of
Canada (NSERC). The second author (RCK) is grateful for the hospitality extended to him
by AMH and her colleagues at Wilfrid Laurier University, by Professors Bill Chen and Arthur Yang at Nankai University and 
by Professor Itaru Terada at the University of Tokyo, and for the financial support making visits to each of these
universities possible.

\vfill
\end{document}